\documentclass[11pt,final,leqno,onefignum,letterpaper]{siamltex}

\usepackage{amsmath}
\usepackage{amssymb}
\usepackage{booktabs}
\usepackage{fixmath}
\usepackage{sansmath}
\usepackage{graphicx}
\usepackage{xcolor}
\usepackage{url}

\title{Dynamics-Adapted Cone Kernels} 

\author{Dimitrios Giannakis\thanks{Center for Atmosphere Ocean Science, Courant Institute of Mathematical Sciences, New York University, New York, New York. 
(dimitris@cims.nyu.edu}).}

\DeclareMathOperator{\Id}{Id}
\DeclareMathOperator{\divr}{div}

\newcommand{\mathbbm}[1]{\text{\usefont{U}{bbm}{m}{n}#1}}
\newtheorem{remark}{{\it Remark}}[section]

\begin{document}
\maketitle
\newcommand{\slugmaster}{%
\slugger{siads}{xxxx}{xx}{x}{x--x}}%slugger should be set to juq, siads, sifin, or siims

\begin{abstract}
We present a family of kernels for analysis of data generated by dynamical systems. These so-called cone kernels feature an explicit dependence on the dynamical vector field operating in the phase-space manifold, estimated empirically through finite-differences of time-ordered data samples. In particular, cone kernels assign strong affinity to pairs of samples whose relative displacement vector lies within a narrow cone aligned with the dynamical vector field. As a result, in a suitable asymptotic limit, the associated diffusion operator generates diffusions along the dynamical flow, and is invariant under a weakly restrictive class of transformations of the data, which includes conformal transformations. Moreover, the corresponding Dirichlet form is governed by the directional derivative of functions along the dynamical vector field. The latter feature is metric-independent. The diffusion eigenfunctions obtained via cone kernels are therefore adapted to the dynamics in that they vary predominantly in directions transverse to the flow. We demonstrate the utility of cone kernels in nonlinear flows on the 2-torus and North Pacific sea surface temperature data generated by a comprehensive climate model.  
\end{abstract}



\pagestyle{myheadings}
\thispagestyle{plain}
\markboth{D.~GIANNAKIS}{DYNAMICS-ADAPTED CONE KERNELS}

\section{Introduction}

Large-scale datasets generated by dynamical systems arise in a diverse range of disciplines in science and engineering, including fluid dynamics \cite{HolmesEtAl96,SreenivasanAntonia97}, materials science \cite{KatsoulakisVlachos03,KatsoulakisEtAl06}, molecular dynamics \cite{DeuflhardEtAl99,NadlerEtAl06}, and geophysics \cite{MajdaWang06,DymnikovFilatov97}. A major challenge in these domains  is to utilize the vast amount of data that is being collected by observational networks or output by large-scale numerical models to understand the operating physical processes, and make inferences about aspects of the system which are not accessible to observation. For instance, in climate atmosphere ocean science (CAOS)  the dynamics takes place in an infinite-dimensional phase space where the coupled nonlinear partial differential equations for fluid flow and thermodynamics are defined, and the observed data correspond to functions of that phase space, such as temperature or circulation measured over a geographical region of interest. There exists a strong need for data analysis algorithms to extract and create reduced representations of the large-scale coherent patterns which are an outcome of these dynamics, including the El Ni\~no Southern Oscillation (ENSO) in the ocean \cite{Trenberth97} and the Madden-Julian Oscillation in the atmosphere \cite{MaddenJulian72}. Advances in the scientific understanding and forecasting capability of these phenomena have potentially high socioeconomic impact. 
  
Despite the high phase space dimension of many systems of interest, their dynamics evolve asymptotically (at long times) on low-dimensional submanifolds of phase space (attractors) \cite{DymnikovFilatov97,KatokHasselblatt97,Arnold98}. It is therefore natural to exploit topological and geometrical aspects of these submanifolds to design data analysis algorithms for data reduction and decomposition, function learning, and other important problems. Here a major challenge stems from the fact that the geometry of the attractive manifolds is nonlinear, meaning that linear variance-optimizing algorithms such as principal components analysis (PCA) \cite{AubryEtAl91,HolmesEtAl96} or singular spectrum analysis (SSA) \cite{BroomheadKing86,VautardGhil89,GhilEtAl02} are likely to perform suboptimally. Indeed, it has been documented in the literature (e.g., \cite{AubryEtAl93,CrommelinMajda04}) that the dynamically significant modes in nonlinear systems are not necessarily those carrying high variance. Moreover, in practical applications involving finite sample counts and short observation intervals, the data manifolds are inherently discrete, and represent only coarse-grained aspects of the underlying geometry at the infinite-sample limit. 

Kernel methods have been extensively used as alternatives to classical linear algorithms to take advantage of nonlinear geometric structures of data. Prominent applications include dimension reduction and feature extraction \cite{ScholkopfEtAl98,BelkinNiyogi03,CoifmanEtAl05,CoifmanLafon06,BelkinNiyogi08,JonesEtAl08}, learning and regularization of scalar or vector-valued functions \cite{BelkinEtAl06,MicchelliPontil05,CaponnettoEtAl08}, and out-of-sample extension of these functions \cite{CoifmanLafon06b}. Here, a common theme is that suitably constructed kernels, i.e., functions measuring a notion of pairwise similarity between data points, lead naturally to diffusion operators which are closely related to the manifold structure of the data. More specifically, because every elliptical diffusion operator induces a Riemannian metric tensor on the data manifold \cite{Rosenberg97,ElworthyEtAl10}, using that operator (or, as is frequently the case, its eigenfunctions) for tasks such as data representation and function learning is tantamount to performing these tasks in a manner compatible with the induced Riemannian geometry.  Applied to datasets with nonlinear manifold structure, such as those generated by complex dynamical systems, this approach has been found to yield more efficient algorithms and meaningful results than linear PCA-type algorithms (e.g.,~\cite{LeeVerleysen07} and references therein).

A major advantage of constructing diffusion operators through kernels is that kernels are defined in the ambient data space, and thus can be used to control the induced Riemannian metric \cite{Berry13} exploiting features of data space which are available for the problem at hand. In the context of dynamical system data, an important feature which is not present in general point clouds is that the samples occur with a time ordering which is the outcome of dynamical flow in phase space. Equivalently, we think of these datasets as being associated with a nowhere-vanishing vector field $ v $ on the attractor, which is intrinsic to the system under study in the sense that it does not depend on how the data manifold is embedded in ambient space. These observations provide a motivation to seek aspects of $ v $ which are empirically accessible in ambient data space, and incorporate them in dynamics-adapted kernels.

Efforts in that direction have been made in a series of papers on so-called nonlinear Laplacian spectral analysis (NLSA) algorithms for decomposition of spatiotemporal data \cite{GiannakisMajda11,GiannakisMajda12a,GiannakisMajda13}, and in independent work by Berry and collaborators \cite{BerryEtAl13,Berry13}. In both of these works, dynamical flow is incorporated in kernels by first embedding the data in a higher-dimensional space (hereafter called embedding space) through Takens' method of delays \cite{PackardEtAl80,Takens81,BroomheadKing86,SauerEtAl91}. Because each point in embedding space corresponds to a segment of dynamical evolution observed over a temporal window, distances in that space depend on the dynamical flow generating the data. This means in turn that the induced Riemannian metric associated with a kernel formulated in embedding space will depend (albeit indirectly) on the dynamical system generating the data. 

In NLSA algorithms, the embedding-space distances are scaled by a factor which is proportional to the distance between temporal nearest neighbors. As will be made precise below, the displacement vector between temporal nearest neighbors is a representation of the dynamical vector field, pushed forward to embedding space. Therefore, in addition to time-lagged embedding, the NLSA kernel also depends directly on the dynamical vector field through its norm. Qualitatively, the result of this dynamics-dependent scaling factor is to assign higher weights to transitory states with large phase space velocity. This feature was found to be particularly useful in extracting dynamically significant modes in systems with metastability \cite{GiannakisMajda12a}.  Moreover, because of the ``non-dimensionalization'' produced by the scaling factor, the NLSA kernel can be naturally extended to process multivariate datasets consisting of components with different physical units \cite{BushukEtAl13}. Yet, in spite of these attractive features, theoretical understanding of the role of the scaling factor in the induced Riemannian metric has so far been lacking.     

Building on the existing work in \cite{GiannakisMajda11,GiannakisMajda12a,GiannakisMajda13}, in this paper we introduce a one-parameter family of kernels which, besides the norm of the dynamical vector field also depend on the angle between that vector field and the displacement vector between the points at which the kernel is evaluated. In particular, the new kernels assign higher weight to pairs of points whose relative displacement lies within a cone with axis parallel to the dynamical vector field. For this reason we refer to the new kernels as cone kernels. The angular influence is controlled by a parameter $ \zeta \in [ 0, 1 ) $, such that the existing NLSA kernel with no angular dependence occurs as the special case $ \zeta = 0 $.  

Here, however, our main focus is the limiting behavior $ \zeta \to 1 $, where the angular influence is maximal. In that limit, the induced Riemannian metric becomes degenerate, and diffusion takes place along the integral curves of the dynamical vector field $ v $ (in the sense of \cite{ElworthyEtAl10}). Because the dynamical integral curves are intrinsic to the system (i.e., do not depend on the ambient-space induced metric), the along-$ v $ property of the $ \zeta \to 1 $ diffusion operator is invariant under arbitrary diffeomorphisms of the data manifold. Moreover, the associated Dirichlet form depends on the directional derivative of functions along $ v $, which is also metric independent. In fact, the only dependence of the limit operator on the ambient-space metric is through a ratio involving its volume form and the norm of $ v $. By virtue of this property, as $ \zeta \to 1 $, the diffusion operators constructed from cone kernels become invariant under a weakly restrictive class of transformations of the data, which includes conformal transformations as a special case. In addition, due to the structure of the Dirichlet form, the associated diffusion eigenfunctions are expected to vary predominantly along directions transverse to $ v $. These eigenfunctions can therefore be employed to carry out tasks such as dimension reduction and feature extraction in a dynamics-adapted manner.

We demonstrate the utility of cone kernels in two applications involving analytically solvable nonlinear flows on the 2-torus and dynamical evolution of sea surface temperature (SST) in the North Pacific ocean in a comprehensive climate model. In the torus application, we explicitly demonstrate the adaptation of the eigenfunctions of the associated diffusion operator to the dynamical flow, and invariance of these eigenfunctions under a non-conformal deformation of the torus. The climate model results suggest that the cone-kernel eigenfunctions produce more efficient (i.e., using fewer basis functions) representation of the temporal variability compared to the existing NLSA kernel, while achieving better timescale separation. 

The plan of this paper is as follows. In section~\ref{secConeKernels}, we lay out the notation, and develop the cone kernel formulation. In section~\ref{secAsymptotics}, we study the behavior of the induced metric and the associated diffusion operator. We present the applications to torus flows and climate model data in sections~\ref{secTorus} and~\ref{secGCM}, respectively, and conclude in section~\ref{secConclusions}. Technical results and proofs are included in an appendix. A movie illustrating the time evolution of SST filtered by diffusion eigenfunctions is provided as supplemental material. A Matlab code used to generate the numerical results in sections~\ref{secTorus} and~\ref{secGCM} is available upon request from the author.   

\section{\label{secConeKernels}Formulation of cone kernels}

\begin{figure}
  \centering
  {\sffamily \sansmath 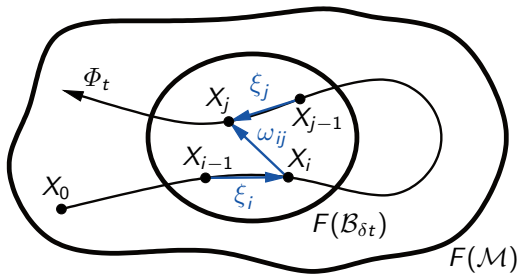}
  \caption{\label{figManifold}Illustration of the $ m $-dimensional data manifold $ \mathcal{ M } $ embedded in $ \mathbb{ R }^n $ through the map $ F $. $ \mathcal{ B }_{\delta t } $ is the open ball used in the asymptotic analysis in Appendix~\ref{appHMetric}.}  
\end{figure}

We consider a scenario where the dynamics is described by a deterministic flow $ \Phi_t : \mathcal{ F } \mapsto \mathcal{ F } $ operating in a phase space $ \mathcal{ F } $, and evolving on a smooth (of class $ C^1 $), compact $ m $-dimensional attractor $ \mathcal{ M } \subseteq \mathcal{ F } $ without boundary.  Moreover, observations are taken uniformly in time with a timestep $ \delta t > 0 $ on the attractor via a $ C^2 $ vector-valued function $ F : \mathcal{ F } \mapsto \mathbb{ R }^n $, leading to a dataset 
\begin{equation}
  \label{eqDiscreteDataset}
  X = \{ X_1, \ldots, X_s \}, \quad \text{with} \quad X_i = F( a_i ), \quad a_i = \Phi_{t_i} a_0, \quad t_i = i \, \delta t, \quad a_0 \in \mathcal{ M }. 
\end{equation}
See Figure~\ref{figManifold} for an illustration. We refer to $ \mathbb{ R }^n $ interchangeably as data space or ambient space, and consider that it is equipped with an inner product $ ( \cdot, \cdot ) $. 

Throughout, we assume that $ \Phi_t $ is sufficiently smooth so that the dynamical vector field $ v $ induced on $ \mathcal{ M } $, defined through
\begin{equation}
  \label{eqV}
  v( f ) = \lim_{t\to 0} ( f( \Phi_t a ) - f( a ) ) / t, \quad \text{with} \quad \text{$ a \in \mathcal{ M } $, $ f \in C^1 \mathcal{ M } $},
\end{equation}
is at least $ C^1$. Moreover, without loss of generality, we assume that $ F $ is one-to-one,  and has full rank on the tangent spaces of $ \mathcal{ M } $; i.e., $ F $ is an embedding of $ \mathcal{ M } $ into $ \mathbb{ R }^n $ (if $ F $ is not an embedding, it is generically possible to construct an embedding through delay-coordinate maps \cite{Takens81,SauerEtAl91}).  Under these conditions, the observation map generates a vector field $ V $ on $ F( \mathcal{ M } ) $ through its derivative $ D F $, viz.\ $ V = DF \, v $.  We denote by $ g $ the Riemannian metric induced on $ \mathcal{ M } $ by pulling back the ambient-space inner product, giving 
\begin{equation}
  \label{eqG}
  g( u, u' ) = ( DF \, u, DF \, u' )
\end{equation}
for any two tangent vectors $ u, u' \in T_a \mathcal{ M } $.

\begin{remark}{\ \rm Even though we do not attempt to extend the analysis presented here to non-smooth manifolds, we note work of Sauer et al.~\cite{SauerEtAl91}, who prove embedding theorems for dynamical systems with fractal attractors. Moreover, in section~\ref{secDiscreteFormulation} we pass to an intrinsically discrete formulation, where the existence of an underlying smooth continuous theory is not required.}   
\end{remark}

Next, consider how linear combinations of $ X_i $ can be used to construct finite-difference (FD) approximations of $ V $. In particular, let $ \delta_p $ be a $ p $-th order FD operator for the first derivative of a function,
\begin{equation}
  \label{eqFD}
  \frac{ df }{ d t } = \frac{ \delta_p f }{ \delta t } + O( \delta t^p ).
\end{equation}
We then have the following Lemma. 
\begin{lemma}\label{lemmaFD}
  The vector in data space $ \xi_i = \delta_p X_i $ corresponds to an $ O( \delta t^p ) $-accurate FD approximation of the pushforward of the dynamical vector field $ v $ evaluated at $ a_i \in \mathcal{ M } $, in the sense that 
  \begin{equation}
    \label{eqFDV}
    DF \, v \rvert_{a_i}  = \xi_i / \delta t + O( \delta t^p ).
  \end{equation}
  Moreover, the data space norm $ \lVert \xi_i \rVert = \sqrt{ ( \xi_i, \xi_i ) } $  provides an $ O( \delta t^p ) $ approximation to the norm of $ v \rvert_{a_i}  $ with respect to the induced metric, namely
  \begin{equation}
    \label{eqFDVNorm}
    \lVert v \rVert_{g,a_i}  := \sqrt{ g( v, v ) } \rvert_{a_i} = \lVert \xi_i  \rVert / \delta t + O( \delta t^p ).
  \end{equation}
\end{lemma}
A proof of this Lemma is included in Appendix~\ref{appFD}. Hereafter, we nominally work with a central FD scheme,
\begin{equation}
  \label{eqFDW}
  \delta_p X_i = \sum_{j=-p}^p w_j X_{i+j}, 
\end{equation}
where $ w_{-p}, \ldots, w_p $ are standard FD weights for central schemes (e.g., \cite{LeVeque07}). However, the asymptotic analysis in section~\ref{secAsymptotics} only depends on the FD accuracy, and also applies, e.g., for backward and forward schemes. 

Recall now that a kernel is a symmetric function $ K : \mathcal{ M } \times \mathcal{ M } \mapsto \mathbb{ R }_+ $ which maps pairs of states in $ \mathcal{ M } $ to a positive number, and, in practical applications, depends only on quantities observed in data space through $ F $. A standard choice in this context is the isotropic Gaussian kernel \cite{BelkinNiyogi03,CoifmanEtAl05,CoifmanLafon06,BelkinNiyogi08},
\begin{equation}
  \label{eqKGauss}
  \bar K_\epsilon( a_i, a_j ) = \exp( - \lVert \omega_{ij} \rVert^2 / \epsilon ),  \quad \text{with} \quad \omega_{ij} = X_j - X_i = F( a_j ) - F( a_i ).
\end{equation}
Here, $ \epsilon $ and $ \omega_{ij} $ are a positive parameter and the displacement vector between data samples $ X_i $ and $ X_j $, respectively. 

Having established the association between $ \xi_i $ and $ v_i $ in Lemma~\ref{lemmaFD}, we seek to modify $ \bar K $ to incorporate information about the dynamical vector field through (i) its norm $ \lVert v \rVert_{g,a_i} $ with respect to the induced metric estimated via~\eqref{eqFDVNorm}; (ii) its angle relative to $ \omega_{ij} $, estimated via
\begin{displaymath}
  \cos \theta_i = \frac{ ( \xi_i, \omega_{ij } ) }{ \lVert \xi_i \rVert \lVert \omega_{ij} \rVert }.
\end{displaymath}
Specifically, introducing a parameter $ \zeta \in [ 0, 1 ) $, we define
\begin{equation}
  \label{eqKCone}
  K_{\delta t, \zeta }( a_i, a_j ) = \exp \left( - \frac{ \lVert \omega_{ij} \rVert^2 }{ \lVert \xi_i \rVert \lVert \xi_j \rVert} [ ( 1 - \zeta \cos^2 \theta_ i )( 1 - \zeta \cos^2 \theta_j ) ]^{1/2} \right).
\end{equation}
For $ \zeta = 0 $, the kernel in~\eqref{eqKCone} reduces to the locally scaled kernel employed in NLSA algorithms \cite{GiannakisMajda12a,GiannakisMajda13}, which features no dependence on the angle between $ \omega_{ij} $ and the phase space velocity vectors, $ \xi_i $ and $ \xi_j $. On the other hand, as $ \zeta $ approaches 1, $ K_{\delta t, \zeta} $ assigns higher affinity to data samples whose relative displacement vector is aligned with either $ \xi_i $ and/or $ \xi_j $. For this reason, we term this two-parameter family of kernels cone kernels. For the remainder of this section we discuss certain properties of cone kernels which will be useful for the asymptotic analysis in section~\ref{secAsymptotics}. 

\subsection{\label{secQualitative}Qualitative features of cone kernels}
 
As is evident from the structure of~\eqref{eqKCone}, the scaling of the pairwise distances $ \lVert \omega_{ij} \rVert $ by $ \lVert \xi_i \rVert \lVert \xi_j \rVert $ results in greater similarity being ascribed to transitory states characterized by large norm of the dynamical vector field. In \cite{GiannakisMajda12a}, this feature was found to be crucial for successful dimensional reduction of a dynamical system with chaotic metastability. In section~\ref{secAsymptotics}, we give a geometrical interpretation of the scaling factor showing that the metric tensor induced on $ \mathcal{ M } $ through $ K_{\delta t, \zeta } $ is invariant under conformal transformations of the data as a result of this scaling. 

In addition, cone kernels with $ \zeta \approx 1 $ provide superior discrimination by assigning greater similarity to those sample pairs whose relative displacement vector is aligned with the dynamical flow. More specifically, given two distinct data samples $ a_i $ and $ a_j $ with $ \omega_{ij} \neq 0 $, the ratio
\begin{displaymath}
  \frac{ K_{\delta t, 1}( a_i, a_j ) }{ K_{\delta t, 0 }( a_i, a_j ) } = \exp\left( \frac{ \lVert \omega_{ij} \rVert^2 (1 - C_{ij}) }{ \lVert \xi_i \rVert \lVert \xi_j \rVert}  \right), \quad C_{ij} = [ ( 1 -  \cos^2 \theta_ i )( 1 - \cos^2 \theta_j ) ]^{1/2} \leq 1
\end{displaymath}
grows exponentially as $ \delta t \to 0 $ whenever $ \cos \theta_i $ and $ \cos \theta_j $ are both not equal to unity. An outcome of this property is that the corresponding diffusion operator $ \upDelta_\zeta $ generates diffusions along  the integral curves of the dynamical vector field $ v $  \cite{ElworthyEtAl10} in the sense that, asymptotically as $ \zeta \to 1 $, $ \upDelta_\zeta f $ vanishes whenever the gradient of $ f $ is parallel to $ v $. Because these curves are intrinsic to the dynamical system generating the data, this property does not depend on the observation function $ F $ and the associated induced metric in~\eqref{eqG}, so long as $ F $ meets the conditions of a manifold embedding \cite{SauerEtAl91}. 

We also mention the utility of cone kernels in situations where $ F $ is a composite map, $ F : \mathcal{ M } \mapsto \mathbb{ R }^{n_1} \oplus \mathbb{ R }^{n_2} $, such that $ F( a ) = ( F_1( a ), F_2( a ) ) $ where both $ F_1 $ and $ F_2 $ are embeddings. This scenario arises in practice when one has access to multivariate observations with distinct physical units, but there is no natural way of choosing a norm for the product space $ \mathbb{ R }^{n_1}  \oplus \mathbb{ R }^{n_2} $. Because the ratio $ \lVert \omega_{ij} \rVert^2 / \lVert \xi_i \rVert \lVert \xi_j \rVert $ is invariant under scaling of the data by a constant (including change of units), cone kernels computed individually for $ F_1 $ and $ F_2 $ can be combined into a single product kernel without having to introduce additional scaling parameters. A climate science application of this technique can be found in \cite{BushukEtAl13}.

\subsection{\label{secLocalBehavior}Local behavior at the basepoint} 

The metric tensor induced at a reference point $ a $ on the data manifold by an exponentially decaying kernel depends strongly on the local rate of decay of the kernel at $ a $. As established by Berry~\cite{Berry13}, in many cases of interest, including the isotropic Gaussian kernel in~\eqref{eqKGauss} and the cone kernel in~\eqref{eqKCone}, that rate is controlled to leading order by the Hessian (second derivative) matrix of the kernel with respect to a coordinate patch covering the reference point, as we now discuss.  

Fixing a basepoint $ a \in \mathcal{ M } $, consider $ K_{\delta t, \zeta }( a, a' ) $ for $ a' $ lying in an exponential neighborhood of $ a $:
\begin{displaymath}
  a' = \exp_a  u, \quad u \in T_a \mathcal{ M }, \quad u = \sum_{\mu=1}^m u^\mu U_\mu.
\end{displaymath}
Here, $ U_1, \ldots U_m $ is a basis of $ T_a \mathcal{ M } $ so that the components $ u^\mu $ are exponential coordinates centered at $ a $. Taking partial derivatives with respect to $ u^\mu $, it is possible to derive the expressions 
\begin{equation}
  \label{eqKDerivatives}
    \left. \frac{ \partial K_{\delta t, \zeta } }{ \partial u^\mu } \right \rvert_{u=0} = 0 \quad \text{and} \quad \left. \frac{ \partial^2 K_{\delta t, \zeta } }{ \partial u^\mu \, \partial u^\nu } \right \rvert_{u=0} = -\frac{ 2 }{ \lVert \xi \rVert^2 } \left( g_{\mu\nu} - \zeta \frac{ \xi^*_\mu \xi^*_\nu }{ \lVert \xi \rVert^2 } \right).
  \end{equation}
In~\eqref{eqKDerivatives}, $ g_{\mu \nu} $ are the components of the ambient-space induced metric $ g $ in~\eqref{eqG} evaluated in the dual basis $ U^{*1}, \ldots, U^{*m} $ with $ U^{*\mu}( U_\nu ) = {\delta^\mu}_\nu $ [see also~\eqref{eqGComponents}]. Moreover, $ \xi^*_\mu $ are the components of the dual vectors $ \xi^* \in T^*_a \mathcal{ M } $ given by pulling back the dual vectors to $ \xi $ with respect to the canonical inner product $ ( \cdot, \cdot ) $ of $ \mathbb{ R }^n $. That is, 
\begin{equation}
  \label{eqXiHat}
  \xi^* := DF^* \, \Xi = \sum_{\mu=1}^m \xi^*_\mu U^{*\mu},
\end{equation}
where $ DF^* : T^*_{F(a)} \mapsto T^*_a \mathcal{ M } $ is the pullback map for dual vectors, and $ \Xi = ( \xi, \cdot ) $. Details of this calculation are provided in Appendix~\ref{appKDerivatives}. 

Equation~\eqref{eqKDerivatives} in conjunction with Lemma~\ref{lemmaFD} leads to an asymptotic expression in the sampling interval $ \delta t $ connecting the Hessian of cone kernels to the dual of the dynamical vector field $ v^* = g( v, \cdot ) = \sum_{\mu=1}^m v^*_\mu U^{*\mu}  $:
\begin{equation}
  \label{eqKHessianV} \left. \frac{ \partial^2 K_{\delta t, \zeta } }{ \partial u^\mu \, \partial u^\nu } \right \rvert_{u=0} = -\frac{ 2 }{ \lVert v \rVert^2_g \, \delta t^2} \left( g_{\mu\nu} - \zeta \frac{ v^*_\mu v^*_\nu }{ \lVert v \rVert^2_g } \right) + O( \delta t^{p-2} ) \quad \text{with} \quad v^*_\mu = \sum_{\nu=1}^m g_{\mu\nu} v^\nu.
\end{equation}
The metric tensor induced on the data manifold by the cone kernel is in fact proportional to the Hessian \cite{Berry13}; in this case $ h_{\mu\nu } $ is given by the negative of the $ O( \delta t^{-2} ) $ coefficient in~\eqref{eqKHessianV}. Below, we study the geometry induced on the data by $ K_{\delta t , \zeta }$ through the associated diffusion operator. 

\begin{remark}
  {\ \rm Equations~\eqref{eqKDerivatives} and~\eqref{eqKHessianV} are unaltered if one replaces the geometric means involving $ \lVert \xi_i \rVert \lVert \xi_j \rVert $ and $ \cos^2\theta_i $ and $ \cos^2 \theta_j $ in~\eqref{eqKCone} with the corresponding arithmetic and harmonic means, i.e., 
    \begin{equation}
      \label{eqKConeAlt}
      K_{\delta t, \zeta }( a_i, a_j ) = \exp \left[ - \frac{ \lVert \omega_{ij} \rVert^2 }{ 4 } \left( \frac{ 1 }{ \lVert \xi_i \rVert^2 } + \frac{ 1 }{ \lVert \xi_j \rVert^2 } \right)  ( 2 - \zeta \cos^2 \theta_ i - \zeta \cos^2 \theta_j ) \right].
    \end{equation}
    For finite $ \delta t $, the behavior of the two kernels will generally differ. In particular, \eqref{eqKCone} is large if either $ \lVert \xi_i \rVert $ or $ \rVert \xi_j \rVert $ are small, whereas \eqref{eqKConeAlt} is large if both $ \lVert \xi_i \rVert $ and $ \rVert \xi_j \rVert $ are small  (similarly for the angular terms). Thus, the kernel in~\eqref{eqKConeAlt} may have higher discriminating power than~\eqref{eqKCone}, but at the same time may require a larger number of samples for stable behavior.}
\end{remark}

\begin{remark}
  {\ \rm In practical applications,  it may be desirable to introduce an additional scaling parameter analogous to $ \epsilon $ in the isotropic Gaussian kernel in~\eqref{eqKGauss}, i.e., 
    \begin{equation}
      \label{eqKConeAlt2}
      K_{\delta t, \epsilon, \zeta }( a_i, a_j ) = \exp \left( - \frac{ \lVert \omega_{ij} \rVert^2 }{ \epsilon \lVert \xi_i \rVert \lVert \xi_j \rVert} [ ( 1 - \zeta \cos^2 \theta_ i )( 1 - \zeta \cos^2 \theta_j ) ]^{1/2} \right).
    \end{equation}
    Apart from an unimportant scaling factor in the Hessian, the presence of $ \epsilon $ in~\eqref{eqKConeAlt2} does not influence the $ \delta t \to 0 $ asymptotics, but provides additional freedom to tune the kernel in situations where one does not have control of the sampling interval.}
\end{remark}

\section{\label{secAsymptotics}The associated diffusion operator and induced metric tensor}

The classical procedure to construct a diffusion operator for geometric analysis of data from a kernel (e.g.,~\cite{BelkinNiyogi03,CoifmanLafon06,BelkinNiyogi08,Berry13}) begins with the introduction of an integral operator acting on scalar functions on the data manifold, which we denote here by $ H_{\delta t, \zeta} $ to make explicit the two parameters appearing in cone kernels. Specifically, 
\begin{equation}
  \label{eqHIntegral}
  H_{\delta t, \zeta} f( a ) = \frac{ 1 }{ \delta t^m } \int_\mathcal{ M } K_{\delta t, \zeta }( a, \cdot ) f  \mu,
\end{equation}
where $ f $ is a sufficiently smooth scalar function on $ \mathcal{ M } $, and $ \mu $ the volume form of the induced metric $ g $ in~\eqref{eqG}. Scaling $ H_{\delta t, \zeta } f $ by the normalization factor 
\begin{equation*}
  \rho_{\delta t, \zeta }( a ) = H_{\delta t, \zeta } \mathbbm{1}( a ) =  \frac{ 1 }{ \delta t^m } \int_\mathcal{ M } K_{\delta t, \zeta }( a, \cdot )  \mu,
\end{equation*}
we obtain the integral operator $ \mathcal{ P }_{\delta t, \zeta}( a ) = H_{\delta t, \zeta} f( a ) / \rho_{\delta t, \zeta}( a ) $. This operator preserves constant functions, i.e., it is an averaging operator. 

Next, let 
\begin{equation}
  \label{eqLContinuous}
  \mathcal{ L }_{\delta t, \zeta } = ( \Id - \mathcal{ P }_{\delta t, \zeta }  ) / \delta t^2.
\end{equation}
This positive-semidefinite operator can be thought of the generator of $ P_{\delta t, \zeta} $. An important property of $ \mathcal{ L }_{\delta t, \zeta } $ is that it annihilates constant functions, i.e., $ \mathcal{ L }_{\delta t, \zeta } \mathbbm{ 1} = 0 $. Following \cite{ElworthyEtAl10}, we refer to such operators as diffusion operators.      

For suitably-defined kernels, $ \mathcal{ L }_{\delta t, \zeta } $ converges to a second-order self-adjoint operator  $ \upDelta_\zeta = \lim_{\delta t \to 0} \mathcal{ L }_{\delta t, \zeta } $. This operator induces a geometry on the dataset in the sense that it corresponds to a unique codifferential operator $ \delta_\zeta $ and Riemannian metric $ h $ such that
\begin{equation}
  \label{eqLapl}
  \Delta_{\zeta} = \delta_\zeta d, \quad \text{where} \quad ( w, df )_h = ( \delta_\zeta w, f )_h 
\end{equation}
for any smooth 1-form field $ w $ and scalar function $ f $. Here, $ ( \cdot, \cdot )_h $ are the canonical (Hodge) inner products for $ p $-form fields associated with $ h $; i.e., 
\begin{equation}
  \label{eqHodge}
  ( w, df )_h := \int_{\mathcal{ M }} h^{-1}( w, df ) \nu, \quad ( \delta_\zeta w, f ) := \int_{\mathcal M } \delta_\zeta( w ) f  \nu,   
\end{equation}
where $ \nu $ and $ h^{-1} $ are the volume form and ``inverse metric'' associated with $ h $, respectively. In particular, we have the following Lemma:
\begin{lemma}
  \label{lemmaHMetric}
  The induced metric tensor $ h $ at $ a \in \mathcal{ M } $ associated with the cone kernels~\eqref{eqKCone} with FD accuracy from Lemma~\ref{lemmaFD} $ p \geq 4 $ is given by 
  \begin{equation}
    \label{eqHMetric}
    h = \frac{ 1 }{ \lVert v \rVert_g^2 } \left( g - \zeta \frac{ v^* \otimes v^* }{ \lVert v \rVert_g^2 } \right),
  \end{equation} 
  where $ g $ is the ambient-space induced metric in~\eqref{eqG}, and $ v^* = g( v, \cdot ) $ the dual dynamical vector field with respect to $ g $. Moreover, the volume forms $ \nu $ and $ \mu $ of $ h $ and $ g $, respectively, are related through the expression
  \begin{equation}
    \label{eqVolH}
    \nu = ( 1 - \zeta )^{1/2} \bar \nu, \quad \text{with} \quad \bar\nu = \mu / \lVert v \rVert_g^m.
  \end{equation}
\end{lemma} 
 A proof of this Lemma can be found in Appendix~\ref{appHMetric}. 
  
 A corollary of Lemma~\ref{lemmaHMetric} is that cone kernels assign a unique metric to equivalence classes of datasets related by conformal transformations. In particular, we say that the datasets associated with the embeddings $ F : \mathcal{ M } \mapsto \mathbb{ R }^n $ and $ \tilde F : \mathcal{ M } \mapsto \mathbb{ R }^{\tilde n } $ are conformally equivalent if there exists a positive function $ r : \mathcal{ M } \mapsto \mathbb{ R }_+ $ such that the induced metric $ \tilde g $ associated with $ \tilde F $ is given by $ \tilde g \rvert_a = r( a ) g \rvert_a $. Because conformally equivalent datasets have the properties $ \lVert v \rVert^2_{\tilde g} = r \lVert v \rVert_g^2 $ and $ \tilde v^* := \tilde g ( v, \cdot ) =  r v^* $, it follows from~\eqref{eqHMetric} that both $ g $ and $ \tilde g $ lead to the same $ h $ metric; an assertion made in section~\ref{secQualitative}. 

\begin{remark}
  {\ \rm \label{remInvariance} The volume form $ \nu $ of $ h $ is invariant under all transformations that preserve the ratio in the right-hand side of~\eqref{eqVolH}. This set of transformations includes conformal transformations as a special case, but also admits more general transformations where changes in Riemannian volume are appropriately compensated by changes in the norm of $ v $. We will return to this point in section~\ref{secZeta0}.}
\end{remark}

A further consequence of~\eqref{eqHMetric} is that $ h $ leads to a contraction of distances in neighborhoods of the data manifold where $ \lVert v \rVert_g $ is uncharacteristically small. Such regions correspond to metastable dynamical regimes separated by rapid transitions with large $ \lVert v \rVert_g $. As remarked in section~\ref{secQualitative}, the ability of locally-scaled kernels (i.e., cone kernels with $ \zeta = 0 $) to discriminate between regimes of this type has been found to be highly beneficial in Galerkin reduced dynamical models with chaotic metastability \cite{GiannakisMajda12a}.

For our purposes, however, of particular interest is the behavior of $ h $ and the associated codifferential and diffusion operators in~\eqref{eqLapl} in the limit $ \zeta \to 1 $, where the directional influence of the dynamical vector field is maximal. In that limit $ h( u, v ) $ vanishes for all tangent vectors $ u \in T_a \mathcal{ M } $, i.e., the induced metric becomes degenerate. Equivalently, the length of the integral curves of $ v $ measured with respect to $ h $ becomes arbitrarily small. In consequence, the diffusion generated by $ \upDelta_1 $ takes place along the integral curves of the dynamical vector field, as we now discuss.

\subsection{\label{secZeta0}Behavior in the $ \zeta \to 1 $ limit}  

Recall that a diffusion operator $ \Delta = \delta d $ acting on $ C^2 $  scalar-valued functions on an $ m $-dimensional manifold $ \mathcal{ M } $ is said to be along a vector field $ v \in T \mathcal{ M } $ if the codifferential $ \delta w $ vanishes for all $ C^1 $ 1-form fields $ w $  lying in the $ ( m - 1 ) $-dimensional subspace $ S_v \subset C^1 T^*\mathcal{ M } $ with $ w( v ) = 0 $ \cite{ElworthyEtAl10}. Intuitively, one thinks of the diffusion process generated by $ \Delta $ to take place along the integral curves of $ v $. A key property of the diffusion operator $ \Delta_\zeta $ in~\eqref{eqLapl} associated with cone kernels is that it is along the dynamical vector field $ v $ asymptotically as $ \zeta \to 1 $, in the sense of the following Lemma.
\begin{lemma}
  \label{lemmaAlong}
  The codifferential operator $ \delta_\zeta $ associated with cone kernels admits the asymptotic expansion 
  \begin{equation}
    \label{eqAlong}
    \delta_\zeta w = \frac{ 1  }{ 1 - \zeta } \bar \delta w + O((1-\zeta)^0), \quad \bar \delta w =  - \divr_{\bar\nu} [ w( v ) v ], 
  \end{equation}
  where $ \divr_{\bar\nu}  $ is the divergence operator associated with the volume form $ \bar\nu $ in~\eqref{eqHMetric}.
\end{lemma}

This Lemma is proved in Appendix~\ref{appAlong}. An asymptotic expansion for the corresponding diffusion operator in~\eqref{eqLapl} follows by setting $ w = df $ in~\eqref{eqAlong} for some $ C^2 $ scalar function $ f $, i.e., 
\begin{equation}
  \label{eqLaplAlong}
  \upDelta_\zeta f = \frac{ 1 }{ 1 - \zeta } \bar \upDelta f + O(( 1- \zeta )^0), \quad \bar \upDelta f = - \divr_{\bar\nu}[ v( f ) v ].
\end{equation}

A consequence of Lemma~\ref{lemmaAlong} is that $ \delta_\zeta w  = O( (1-\zeta)^0 ) $ if $ w \in S_v $, but $ \delta_\zeta w'  = O( ( 1- \zeta )^{-1} ) $ if $ w' $ lies outside of $ S_v $. As a result,  the norm ratio 
\begin{equation}
  \label{eqAsymptoticAlong}
  \frac{ \lVert \delta_\zeta w \rVert_h }{ \lVert \delta_\zeta w' \rVert_h } = \frac{ ( \delta_\zeta w, \delta_\zeta w )^{1/2} }{ ( \delta_\zeta w', \delta_\zeta w' )^{1/2} } = O(( 1- \zeta ))
\end{equation}
tends to zero as $ \zeta \to 1 $ for all nonzero square-integrable $ 1 $-form fields $ w' \notin S_v $ for which $ \bar \delta w' $ is nonzero.  We interpret~\eqref{eqAsymptoticAlong} as an asymptotic along-$ v $ property of $ \upDelta_\zeta $. 

\begin{remark}
  {\ \rm In general,  $ \upDelta_\zeta f = \delta_\zeta d f = - \divr_\nu \grad_h f $ depends on the metric through the volume form $ \nu $, as well as explicitly through the gradient $ \grad_h f = h^{-1}( df, \cdot ) $. In the $ \zeta \to 1 $ limit, the latter is replaced by the directional derivative $ v( f ) $ in~\eqref{eqLaplAlong}, which is metric-independent.  Thus, the only metric dependence of the limit operator $ \bar \upDelta $ is through the volume form $ \bar \nu $. According to Remark~\ref{remInvariance}, the latter is invariant under a set of transformations of the data which includes conformal transformations as a subset.} 
\end{remark}  

\subsection{\label{secDiscreteFormulation}Discrete formulation}

Discrete analogs of the diffusion operator~\eqref{eqLapl} arise naturally in the framework of discrete exterior calculus (DEC; e.g., \cite{DesbrunEtAl05,HeinEtAl05,ZhouBurges08,GradyPolimeni10}). In this setting, the spaces of scalar functions and 1-form fields appearing in~\eqref{eqLapl} are replaced by functions $ f( a_i ) = f_i $ and $ w( [a_i a_j ] ) = w_{ij} $ defined on the vertices $ a_i $ and edges $ [ a_i a_j ] $, respectively, of a graph formed by the $ s $ sampled states in~\eqref{eqDiscreteDataset}. These function spaces are equipped with weighted inner products, 
\begin{equation}
  \label{eqDiscreteHodge}
  ( f, f' )_P = \sum_{i=1}^s \pi_i f_i f'_i, \quad ( w, w' )_P = \sum_{i,j=1}^s \pi_i P_{ij} w_{ij} w'_{ij} / 2, 
\end{equation}
which are the discrete counterparts of~\eqref{eqHodge}. Also, a difference operator $ \hat d $ is is introduced mapping vertex to edge functions via $ \hat d f([a_i a_j]) = f_j-f_i$. The associated discrete codifferential $ \hat \delta_{\delta t, \zeta } $ and diffusion operator $ L_{\delta t, \zeta } $ (which in the context of cone kernels depend on both $ \delta t $ and $ \zeta $) are then defined in direct analogy with~\eqref{eqLapl}:
\begin{equation}
  \label{eqDiscreteLapl}
  L_{\delta t, \zeta} = \hat \delta_{\delta t, \zeta} \hat d, \quad ( w, \hat df )_P = ( \hat \delta_{\delta t, \zeta} w, f )_P.
\end{equation}
  
Even though the inner products and associated diffusion operator in~\eqref{eqDiscreteLapl} exist independently of a continuous theory, here we seek to construct $ L_{\delta t, \zeta } $ such that, asymptotically as $ s \to \infty $, it inherits the conformal invariance and along-$ v $ properties  established in Lemmas~\ref{lemmaHMetric} and~\ref{lemmaAlong}. To that end, we employ the diffusion map (DM) algorithm of Coifman and Lafon \cite{CoifmanLafon06}. In DM, the inner product weights $ P_{ij} $ in~\eqref{eqDiscreteHodge} are the elements of a Markov matrix whose state space is the discrete dataset in~\eqref{eqDiscreteDataset}, constructed via the sequence of operations 
\begin{equation}
  \label{eqDM}
  K_{ij} = K_{\delta t, \zeta}( a_i, a_j ), \quad \tilde K_{ij} = \frac{ K_{ij} }{ \left( \sum_{k=1}^s K_{ik} \right)^\alpha \left( \sum_{k=1}^s K_{jk} \right)^\alpha}, \quad P_{ij} = \frac{ \tilde K_{ij} }{ D_{ii} D_{jj} }.
\end{equation}
Here, $ \alpha $ is a real parameter, and $ D $ a diagonal degree matrix with $ D_{ii} = \sum_{k=1}^s \tilde K_{ik} $. The inner-product weights $ \pi_i $  are given by the invariant distribution of that Markov matrix, i.e., $ \sum_{i=1}^s \pi_i P_{ij} = \pi_j $. With these definitions, it follows that
\begin{equation}
  \label{eqLaplDiscrete} \hat \delta_{\delta t, \zeta } w( a_i ) = \sum_{j=1}^s P_{ij} (w([a_j a_i]) - w([a_i a_j ] ) / 2, \quad L_{\delta t, \zeta } = I - P.
\end{equation}
  
The normalization step to obtain $ \tilde K_{ij} $ in~\eqref{eqDM}, which is not present in standard graph Laplacian algorithms (e.g., \cite{BelkinNiyogi08}), controls the influence of the sampling density with respect to the volume form $ \nu $ in approximations of the manifold integral in~\eqref{eqHIntegral} by discrete sums of the form $ \sum_{j=1}^s \tilde K_{ij} f( a_j ) $.  In particular, under relatively weak assumptions on $ \mathcal{ M } $ and the embedding map $ F $, it can be shown \cite{CoifmanLafon06} that for $ \alpha = 1 $, $ \delta t \to 0 $, and $ s^{1/2} \delta t^{2+m/2} \to \infty $ \cite{Singer06}, the discrete diffusion operator  converges pointwise to $ \upDelta_\zeta $, in the sense that $ \lvert \upDelta_\zeta f( a ) - L_{\delta t, \zeta } f( a ) / \delta t^2 \rvert \to 0 $ for all $ a \in \mathcal{ M } $ and sufficiently smooth $ f $. This feature is particularly desirable for our purposes, for the sampling density is dictated by the dynamical flow $ \Phi_t $ and observation function $ F $, and we are interested in targeting $ \mathcal{ L }_{\delta t, \zeta } $ without using a priori information about $ \Phi_t $ and $ F $.  We therefore adopt $ \alpha = 1 $ DM in all of the experiments of sections~\ref{secTorus} and~\ref{secGCM} ahead. 

\subsection{\label{secEigenfunctions}Diffusion eigenfunctions}

As mentioned in the introduction, diffusion operators such as $ L_{\delta t, \zeta } $ in~\eqref{eqLaplDiscrete} are useful for a wide range of data analysis tasks, including dimension reduction, feature extraction, and regularization \cite{ScholkopfEtAl98,BelkinNiyogi03,CoifmanEtAl05,CoifmanLafon06,BelkinNiyogi08,JonesEtAl08,BelkinEtAl06,MicchelliPontil05,CaponnettoEtAl08,CoifmanLafon06b}. Below, we focus on a particular aspect of $ L_{\delta t, \zeta } $, namely its eigenfunctions $ \phi_i $, defined through
\begin{equation}
  \label{eqPhi}
  L_{\delta t, \zeta } \phi_i = \lambda_i \phi_i, \quad \phi_i = ( \phi_{1i}, \ldots, \phi_{si} )^T, \quad 0 = \lambda_0 < \lambda_1 \leq \lambda_2 \leq \cdots. 
\end{equation}
Diffusion eigenfunctions are traditionally used to create low-dimensional parameterizations of data of the form $ a \in \mathcal{ M } \mapsto \phi( a ) = ( \phi_{i_1}( a ), \ldots, \phi_{i_l}( a ) ) \in \mathbb{ R }^l $, with rigorous embedding results established in the continuous limit \cite{BerardEtAl94,JonesEtAl08,Portegies13}. A somewhat different perspective, adopted in NLSA algorithms \cite{GiannakisMajda11,GiannakisMajda12a,GiannakisMajda13}, is to associate low-dimensional subspaces spanned by the leading $ \phi_i $ with spaces of temporal patterns through the time series 
\begin{equation}
  \label{eqPhiTS} \tilde \phi_i( t_j ) = \phi_i( a_j ) = \phi_{ji}, 
\end{equation}
and (in the spirit of \cite{AubryEtAl91}) extract spatiotemporal modes of variability through singular value decomposition of the data projected onto those eigenfunctions. Clearly, in both approaches, the relationship of the basis to the underlying dynamical flow plays a major role on algorithm performance. 

To gain insight on the influence of the along-$ v $ property of diffusion operators associated with cone kernels on their eigenfunctions, it is useful to consider the Dirichlet form associated with $ \upDelta_\zeta $,
\begin{displaymath}
  \mathcal{ E }_\zeta( f ) = ( f, \upDelta_\zeta f )_h.
\end{displaymath}
In particular, it follows from the asymptotic expansion in~\eqref{eqAlong}, in conjunction with the divergence theorem and the assumption in section~\ref{secConeKernels} that $ \mathcal{ M } $ has no boundary, that
\begin{equation}
  \label{eqDirichlet}
  \mathcal{ E }_\zeta( f ) = ( 1 - \zeta )^{-1/2 } \bar{\mathcal{ E }}( f ) + O( ( 1 - \zeta )^{3/2} ), \quad \text{with} \quad \bar{\mathcal{E}}( f ) = \int_\mathcal{ M } [ v( f ) ]^2 \bar \nu.
\end{equation} 
Therefore, as $ \zeta \to 1 $, $ \mathcal{ E }_\zeta $ assigns low energy to functions which (on average) have large directional derivative $ v( f ) $ along the dynamical flow. Because  $ \mathcal{ E }_\zeta( \phi_i ) $ is equal to the corresponding eigenvalue $ \lambda_i $ for normalized eigenfunctions with  $ \lVert \phi_i \rVert_h = 1 $, we expect the leading (small-$ \lambda_i $) diffusion eigenfunctions to vary predominantly in directions transverse to the integral curves of $ v $.  This property implies strong adaptivity of $ \phi_i $ to the dynamical flow, and is also independent of the embedding map $ F $ and the associated induced metric in~\eqref{eqG}. Below, we demonstrate these properties in numerical experiments.

\section{\label{secTorus}Dynamical systems on the 2-torus}

To explicitly illustrate the key features of cone kernels, we begin with a low-dimensional application where the phase space manifold $ \mathcal{ M } $ is the 2-torus.  Denoting by $ ( \theta^1, \theta^2 ) $ the azimuthal and polar angles on the 2-torus, respectively, we consider a two-parameter family of dynamical vector fields 
\begin{equation}
  \label{eqTorusDynamical}
  v = \sum_{\mu=1}^2 v^\mu \frac{ \partial }{ \partial \theta^\mu }, \quad \text{with} \quad v^1 = 1 + ( 1 - \beta )^{1/2} \cos \theta^1, \quad   v^2  = \Omega ( 1 - ( 1 - \beta )^{1/2} \sin \theta^2 ).
\end{equation}
Here, $ \Omega  $ is a positive frequency parameter which is set to an irrational number to produce a dense cover of the torus. Moreover, the parameter $ \beta \in ( 0, 1 ] $ controls the nonlinearity of the flow. Specifically,  $ \beta = 1 $ corresponds to a linear flow, but when $ \beta < 1 $ the flow ``slows down'' at $ (\theta^1, \theta^2 ) \sim ( \pi, \pi/2 ) $ and ``speeds up'' at $ ( \theta^1, \theta^2 ) \sim ( 0, -\pi/2 ) $. The orbit $ ( \dot \theta^1, \dot \theta^2 ) = ( v^1, v^2 ) $ passing through $ ( \theta^1, \theta^2  ) = ( 0, 0 ) $ at time $ t = 0 $ is given by
\begin{displaymath}
  \tan ( \theta^1 / 2 ) = [ 1 + ( 1 - \beta )^{1/2} ] \beta^{-1/2}  \tan( \beta t / 2 ), \quad \cot( \theta^2 / 2 ) = ( 1 - \beta )^{1/2} + \beta^{1/2} \cot( \beta^{1/2} t / 2 ).
\end{displaymath} 
 
\begin{figure}
  \centering
  \includegraphics{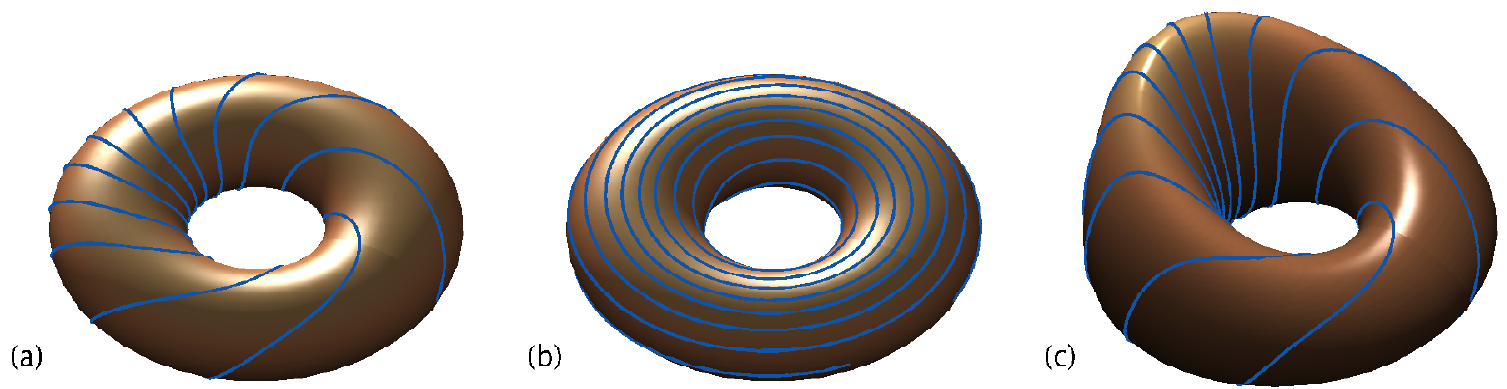}
  \caption{\label{figTrajectories}Sample trajectories for the dynamical systems~\eqref{eqTorusDynamical} on 2-tori. (a) Model~I: $ \Omega = 30^{1/2} $, no deformation. (b) $ Model~II: \Omega = 30^{-1/2} $, no deformation. (c) Model~I$'$: $ \Omega = 30^{1/2} $, non-conformal deformation with $ \gamma = 0.3 $.}
\end{figure}

Setting $ \beta = 0.5 $ throughout, we consider the cases $ \Omega = 30^{1/2} $ and $ 30^{-1/2}  $, referred to here as Models~I and~II, respectively. The resulting trajectories in data space corresponding to the standard embedding of the 2-torus, 
\begin{gather*}
  F: \mathcal{ M } \mapsto \mathbb{ R }^3, \quad F( a ) = ( x^1, x^2, x^3 ), \\
  x^1 = ( 1 + R \cos \theta^2( a ) ) \cos \theta^1( a ), \quad x^2 = ( 1 + R \cos \theta^2( a ) ) \sin\theta^1( a ), \quad x^3 = \sin \theta^2( a ), 
\end{gather*}
with $ R = 1 /2 $ are illustrated in Figures~\ref{figTrajectories}(a) and~\ref{figTrajectories}(b), respectively. In addition to the standard embedding, we also study a non-conformally deformed embedding $ F_\gamma(  a ) = ( y^1, y^2, y^3 ) $ of the 2-torus into $ \mathbb{ R }^3 $, where the $ y^3 $ coordinate is a stretched version of $ x^3 $ with a non-uniform scaling factor, i.e., 
\begin{displaymath}
  y^1 = x^1, \quad \tilde y^2 = x^2, \quad y^3 = x^3 e^{\gamma z}, \quad z = ( 1 + R - x^1 )( 1 + x^3 ).
\end{displaymath}
The deformed system with $ \Omega = 30^{1/2} $ and $ \gamma = 0.3 $, which we refer to as Model~I$^\prime$, is shown in Figure~\ref{figTrajectories}(c). In all cases, we equip the ambient data space  with the canonical Euclidean inner product.

For each of the models in Figure~\ref{figTrajectories}, we generated datasets with $ s = \text{64,000} $ samples taken at a sampling interval $ \delta t = 2 \pi / ( S \min\{ 1, \Omega \} ) $ where $ S = 500 $ is an integer parameter controlling the number of samples in each quasi-period. We computed the FD approximations of $ v $ using the central scheme in~\eqref{eqFDW} with fourth-order accuracy as required by Lemma~\ref{lemmaFD} [i.e., the FD weights are $(w_{-2},\ldots,w_2) = ( 1/12, -2/3, 0, 2/3, -1/12)$]. With this FD scheme, the phase space velocity norm ratio $ \max \{ \lVert \xi_i \rVert \} / \min \{ \lVert \xi_i \rVert \} $ for Models~I, II, and I$^\prime$ was of order 20, 35, and 40, respectively. Thus, the influence of the $ \lVert \xi_i \rVert \lVert \xi_j \rVert $ scaling factors in the cone kernels~\eqref{eqKCone} is expected to be significant. 

We evaluated the discrete diffusion operator $ L_{\delta t, \zeta } $ from cone kernels with $ \zeta = 0 $ and 0.995 (corresponding to no influence and strong influence of the directionality of $ \xi $, respectively) using the DM procedure in~\eqref{eqDM} with $ \alpha = 1 $. For comparison, we also computed the diffusion operator associated with the isotropic Gaussian kernel in~\eqref{eqKGauss} with $ \epsilon = 0.1 $. To limit memory usage, in all cases we truncated the pairwise kernel evaluations to $ b = 2000 $ nearest  neighbors [in the sense of $ K_{ij} $ in~\eqref{eqDM}] per sample. The kernel values at truncation were no greater than $ O( 10^{-7} ) $, indicating that truncation has negligible impact on the numerical results.    

\begin{figure}
  \centering
  \includegraphics{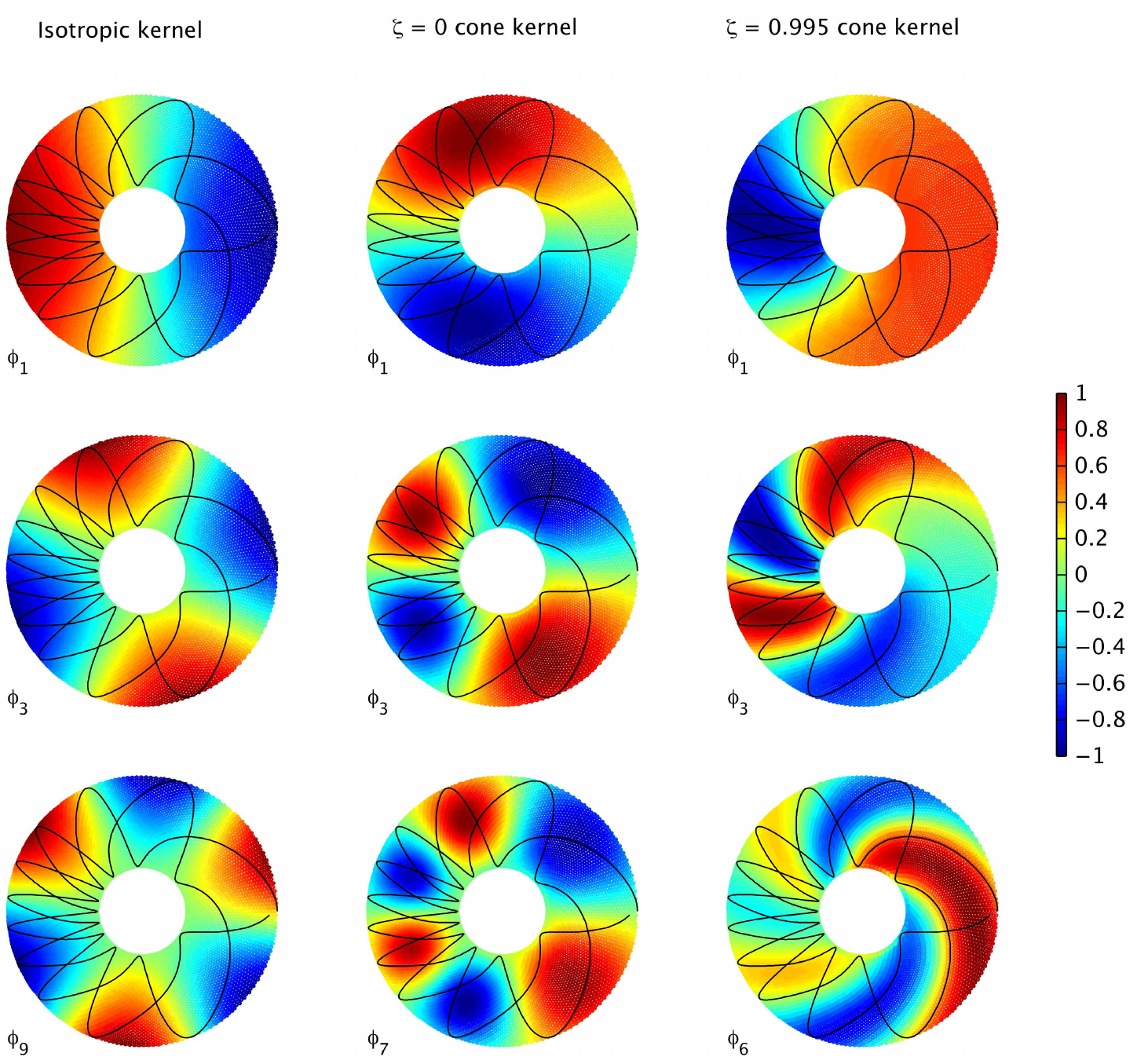}
  \caption{\label{figPhiI}Scatterplots of diffusion eigenfunctions for the dynamical system on the 2-torus with $ \Omega = 30^{1/2} $ obtained using DM with $ \alpha = 1 $ in conjunction with the isotropic Gaussian kernel~\eqref{eqKGauss}, and the cone kernels~\eqref{eqKCone} with $ \zeta = 0 $ and 0.995. A portion of the dynamical system trajectory is plotted in a black line for reference.}  
\end{figure}

\begin{figure}
  \centering
  \includegraphics{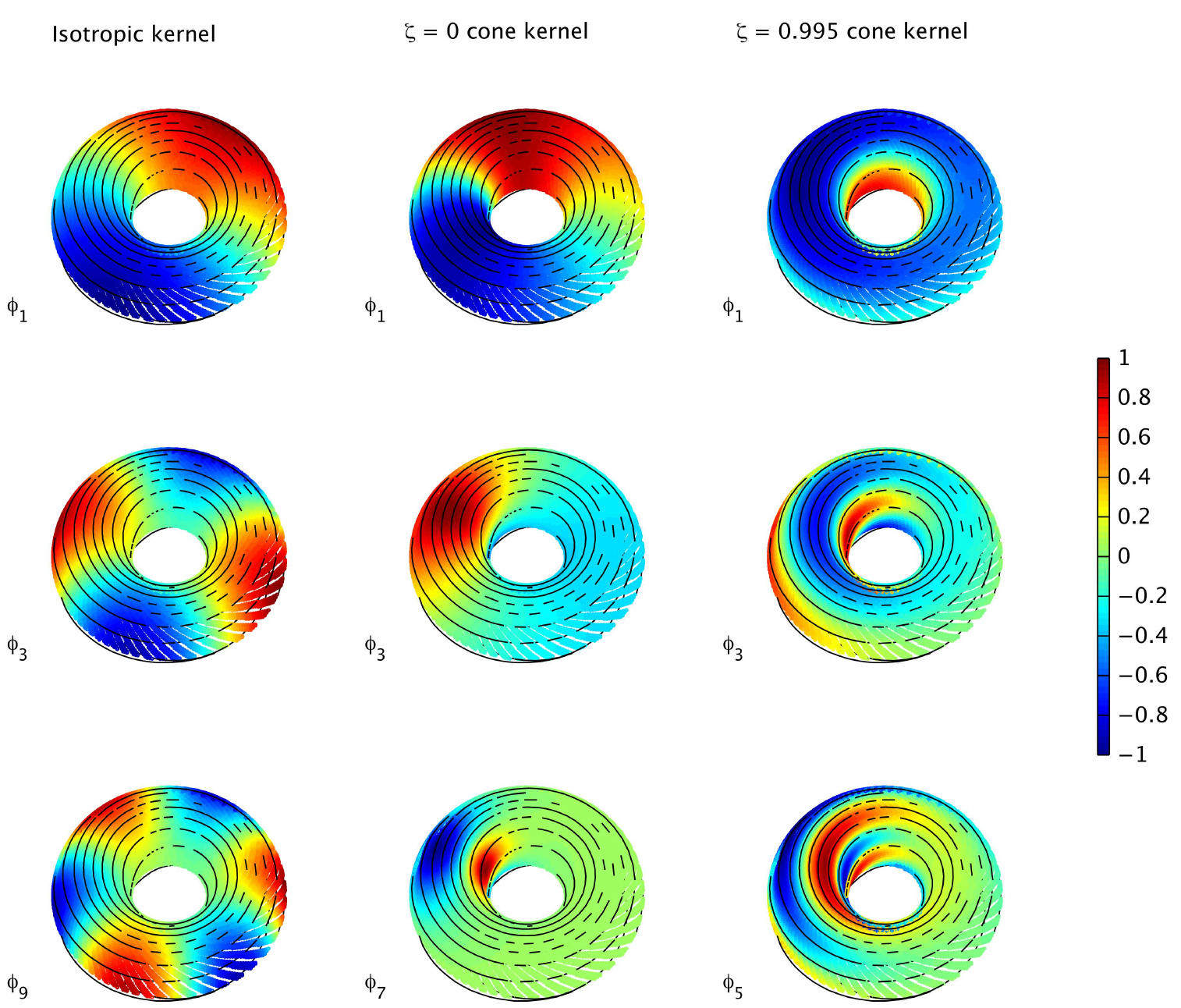}
  \caption{\label{figPhiII}Same as Figure~\ref{figPhiI} but for the dynamical system with $ \Omega = 30^{-1/2} $.}  
\end{figure}

\begin{figure}
  \centering
  \includegraphics{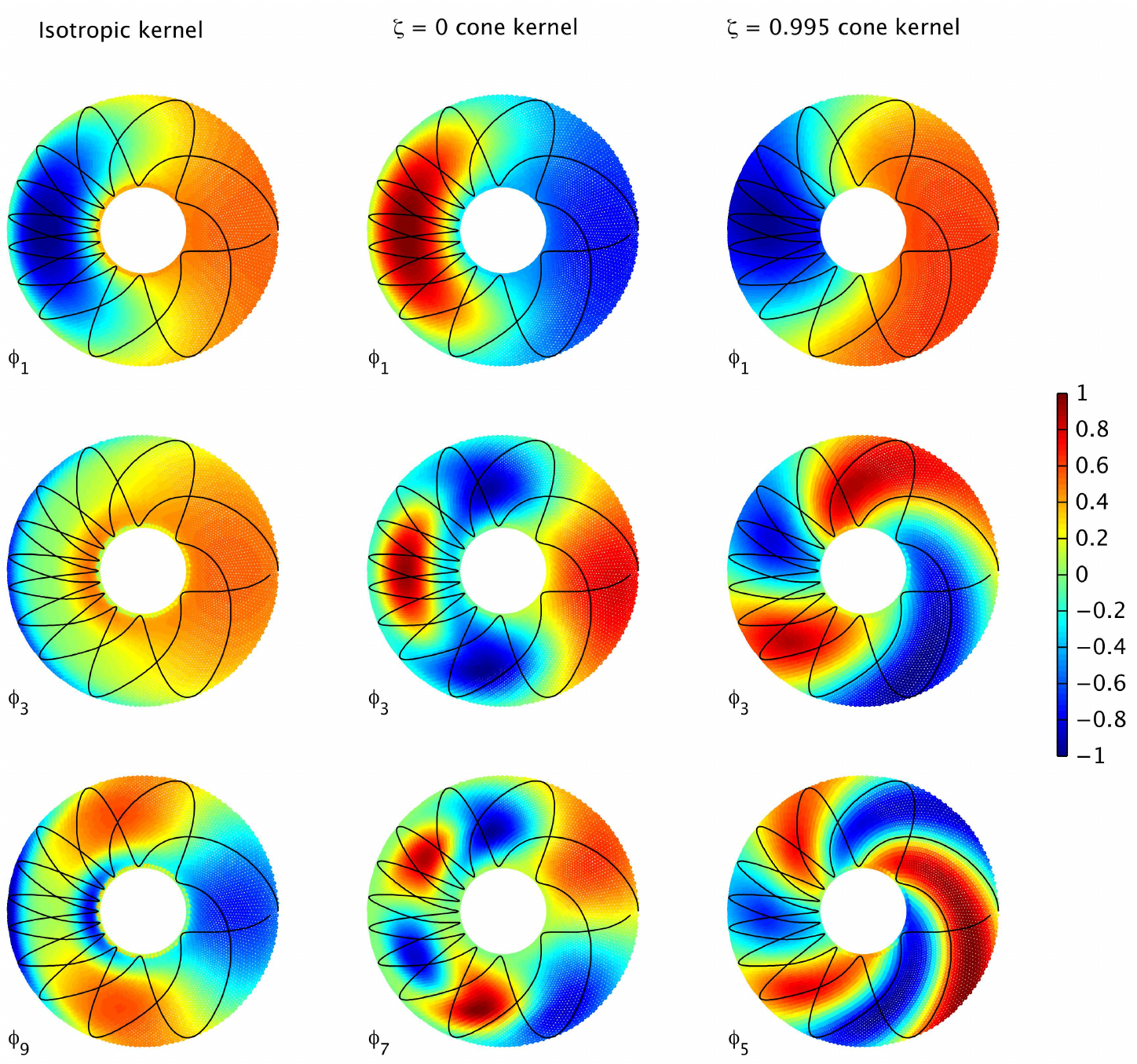}
  \caption{\label{figPhiIDeformed}Same as Figure~\ref{figPhiI} but for the deformed torus with deformation parameter $ \gamma = 0.3 $ [see Figure~\ref{figTrajectories}(c)].}  
\end{figure}

Representative eigenfunctions for Models I, II, and I' are shown in Figures~\ref{figPhiI}--\ref{figPhiIDeformed}. To test for convergence of our results to the continuous limit, we performed a series of long runs, where $ s $ was increased eightfold, $ \delta t$ decreased twofold, and $ b $ set to 10,000 (in the isotropic Gaussian kernel case $ \epsilon $ was decreased fourfold). The diffusion eigenfunctions from the long runs (not shown here) were in good agreement with those displayed in Figures~\ref{figPhiI}--\ref{figPhiIDeformed}. We also performed eigenfunction calculations using a first-order backward FD scheme so that $ \xi_i = X_i - X_{i-1} $, and found again only minor changes relative to the results with fourth-order FD accuracy.  It therefore appears that in this setting with dense sampling the accuracy required by Lemma~\ref{lemmaHMetric} is not crucial to recover the salient features of the eigenfunctions.

\subsection{Properties of the diffusion eigenfunctions}

First, consider the eigenfunctions for Models~I and~II. According to~\cite{CoifmanLafon06}, the results of DM with the isotropic Gaussian kernel should only depend on the embedding $ F $, and should be independent of the sampling density on $ \mathcal{ M } $ induced by the dynamical flow. Indeed, as expected from theory, the isotropic-kernel eigenfunctions in Figures~\ref{figPhiI} and~\ref{figPhiII} are essentially equivalent, despite the fact that the underlying dynamical system trajectories are qualitatively different. 

On the other hand, because cone kernels depend explicitly on the dynamical flow through $ \xi $ [see~\eqref{eqKCone}], the eigenfunctions for Models~I and~II differ. In the $ \zeta = 0 $ results for Model~I shown in Figure~\ref{figPhiI}, the eigenfunction wavecrests are compressed in the portion of the torus where the flow evolves slowly, and rarefied in the transitory region characterized by large $ \lVert v \rVert_g $ (compare, e.g., the $ \phi_3 $ and $ \phi_5 $ eigenfunctions obtained via the isotropic kernel and the $ \zeta = 0 $ cone kernel). As a result, these eigenfunctions have higher discriminating power in the regions of $ \mathcal{ M } $ where the system evolves more slowly. The Model~II eigenfunctions exhibit a qualitatively similar behavior.   

The along-$ v $ property of the diffusion operators with $ \zeta \approx 1 $, and the resulting adaptation of the eigenfunctions to the dynamical flow expected on the basis of~\eqref{eqDirichlet}, is manifested in the right-hand columns of Figures~\ref{figPhiI} and~\ref{figPhiII}. There, the eigenfunctions vary predominantly in directions transverse to the dynamical flow, resulting in characteristic swirl patterns for Model~I and azimuthal streak-like patterns for Model~II. 

Because the integral curves of $ v $ are metric independent, this property is robust against changes in the embedding function $ F $.  Indeed, in the Model~I' results in Figure~\ref{figPhiIDeformed}, it is only the $ \zeta = 0.995 $ eigenfunctions which remain qualitatively similar to those in Figure~\ref{figPhiI}. There, the dataset deformation has left a clear imprint on the eigenfunctions. On the other hand, the $ \zeta = 0.995 $ eigenfunctions retain the dynamics-adaptation featuring the characteristic swirl patterns following the dynamical flow. The structure of the eigenfunctions transversely to the flow does exhibit some changes relative to Figure~\ref{figPhiI}, but these are significantly weaker compared to the isotropic Gaussian kernel and $ \zeta = 0 $  cone kernel examples.   
  
\section{\label{secGCM}North Pacific SST data from a comprehensive climate model}

In this experiment, we apply cone kernels to SST data in the North Pacific sector of the Community Climate System Model version 3 \cite{CollinsEtAl06}. The dataset, which was studied in \cite{GiannakisMajda12b} via NLSA, consists of monthly-averaged samples of SST in the rectangular domain 120$^\circ$E--110$^\circ$W and 20$^\circ$N--65$^\circ$N spanning a 900-year interval.\footnote{The dataset is available at the Earth System Grid repository, \url{http://www.earthsystemgrid.org}, where it is designated CCSM3 integration b30.004.} We work throughout with the model's nonuniform native ocean grid of $ 1^\circ $ nominal horizontal resolution. The number of gridpoints in the analysis domain is $ d = 6671 $. 

A major component of North Pacific SST variability is due to the annually varying solar forcing (the seasonal cycle). The latter is superposed to low-frequency (interannual to decadal) variability patterns, the most prominent of which are the Pacific Decadal Oscillation (PDO) \cite{MantuaHare02} and the North Pacific Gyre Oscillation (NPGO) \cite{DiLorenzoEtAl08}. A signature of ENSO (which is most prominent in the tropical Pacific \cite{Trenberth97}) is also present in this domain. We refer the reader to \cite{GiannakisMajda12b} for further details on these modes of variability extracted via NLSA. Here, our objective is to compare the diffusion eigenfunctions of cone kernels with $ \zeta \approx 1 $ to those with $ \zeta = 0 $, which are equivalent to the earlier NLSA kernels (aside from the fourth-order accurate FD scheme used here in place of the first-order backward scheme in \cite{GiannakisMajda12b}). 

Following \cite{GiannakisMajda12b,GiannakisMajda13}, to ``Markovianize'' the time-dependent solar forcing in the data, and induce timescale separation in the diffusion eigenfunctions, we first embed the data to a higher-dimensional space via Takens' method of delays \cite{PackardEtAl80,Takens81,BroomheadKing86,SauerEtAl91}. That is, we map each spatial snapshot $ x_i \in \mathbb{ R }^d $ to a spatiotemporal sequence $ X_i = ( x_i, x_{i-1}, \ldots, x_{i-(q-1)} ) \in \mathbb{ R }^n $, where $ q $ is an integer parameter measuring the embedding window length (in months). Here, we work with a two-year embedding window, $ q = 24 $, which corresponds to a data space dimension $ n = q d = \text{160,104} $. Qualitatively similar results can be obtained with embedding windows in the interval 1--5 years. Similarly to the torus examples of section~\ref{secTorus}, we evaluate $ \xi $ using a central fourth-order FD scheme. The number of samples $ s $ available for analysis after embedding and removal of two samples in the beginning and end of the simulation interval (to compute $ \xi_1 $ and $ \xi_s$ ) is 10,773. Unlike the torus experiments, $ s $ is small-enough in this case so as not to require nearest-neighbor truncation. We employ an area-weighted data space inner product given by $ ( X_i, X_j ) = \sum_{\nu=1}^n A_{\nu \bmod d }  X_i^\nu X_j^\nu $, where $ A_1, \ldots, A_d $ are the grid cell areas in the analysis domain. Prior to analysis we center the data by subtracting the temporal mean $ \bar x = \sum_i x_i $ from each snapshot $ x_i $ to produce temperature anomaly snapshots. Note that centering is done mainly for visualization purposes, and does not influence the kernel values in~\eqref{eqKGauss} and~\eqref{eqKCone}.   

\begin{figure}
  \centering\includegraphics[scale=.92]{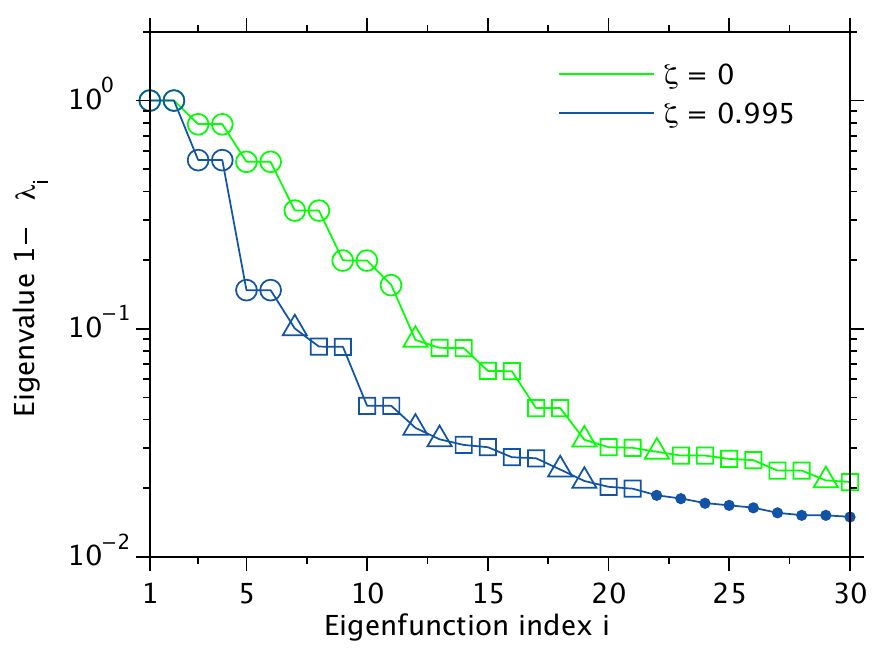}
  \caption{\label{figCCSMLambda}Eigenvalues $ \lambda_i $ for CCSM North Pacific SST data obtained via cone kernels with $ \zeta = 0 $ and $ 0.995 $, normalized so that $ 1 - \lambda_1 = 1 $. To highlight the differences between the two spectra, $ 1 - \lambda_i $ values are plotted in a logarithmic scale. Periodic, low-frequency, and intermittent eigenfunctions are indicated using $\bigcirc$,  $ \bigtriangleup $, and $ \Box $ markers, respectively. Solid markers correspond to eigenfunctions whose temporal character does not belong in these families.}
\end{figure}

\begin{figure}
  \centering\includegraphics[scale=.9]{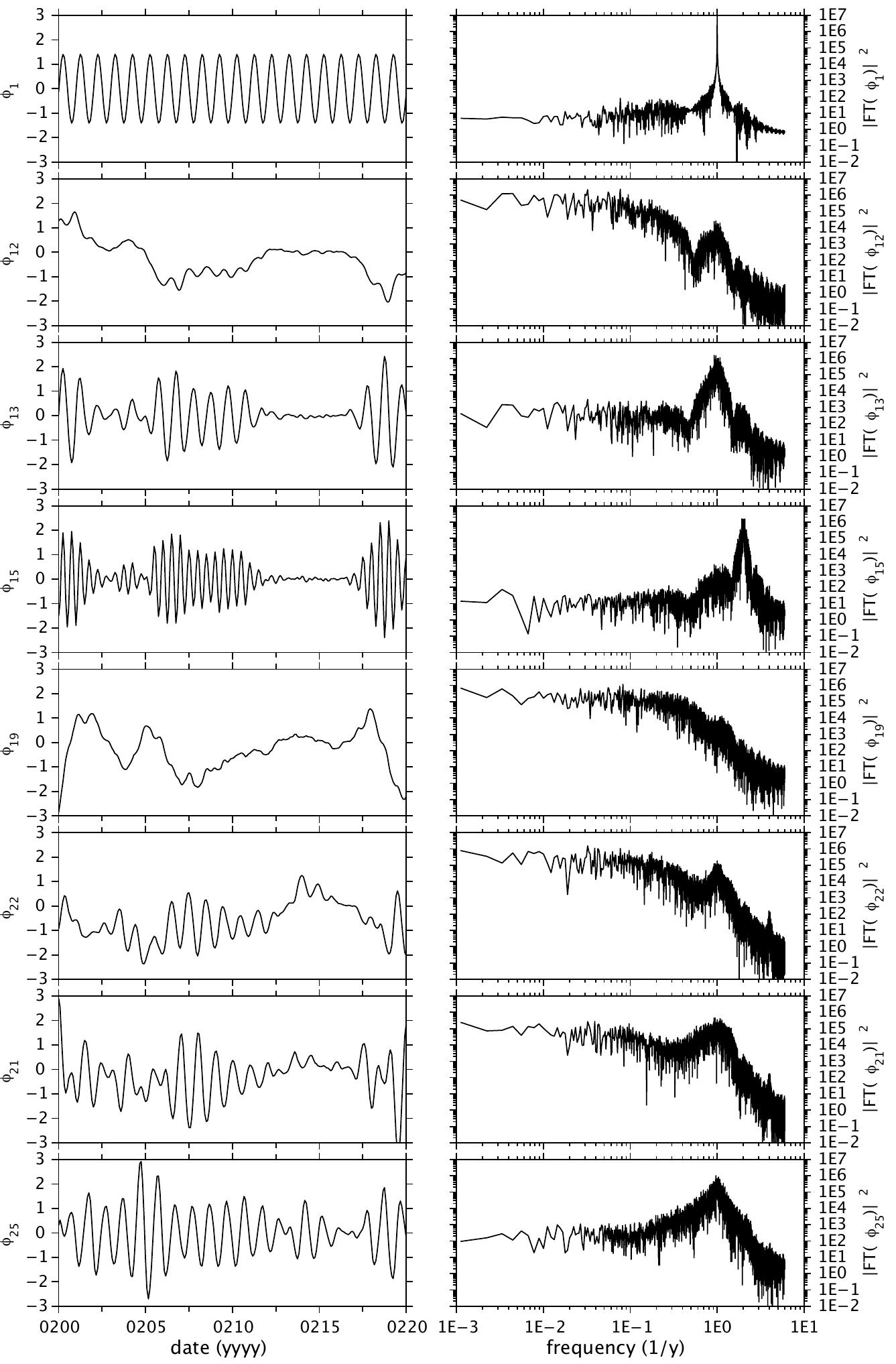}
  \caption{\label{figCCSMPhiAT}Representative eigenfunctions for CCSM North Pacific SST data obtained via the cone kernel with $ \zeta = 0 $. The left-hand panels show 20-year portions of the time series $ \tilde \phi_i( t ) = \phi_i( a_t ) $. The right-hand panels are frequency spectra computed via the discrete Fourier transform of $ \tilde \phi_i $. The eigenfunctions in this Figure have been chosen so as to qualitatively match as close as possible those of Figure~\ref{figCCSMPhiCone}. Note the timescale mixing in $ \{ \phi_{21}, \phi_{22} \} $ compared to the corresponding $ \zeta = 0.995 $ eigenfunctions, $ \{ \phi_{13}, \phi_{14} \} $, in Figure~\ref{figCCSMPhiCone}.}
\end{figure}

\begin{figure}
  \centering\includegraphics[scale=.9]{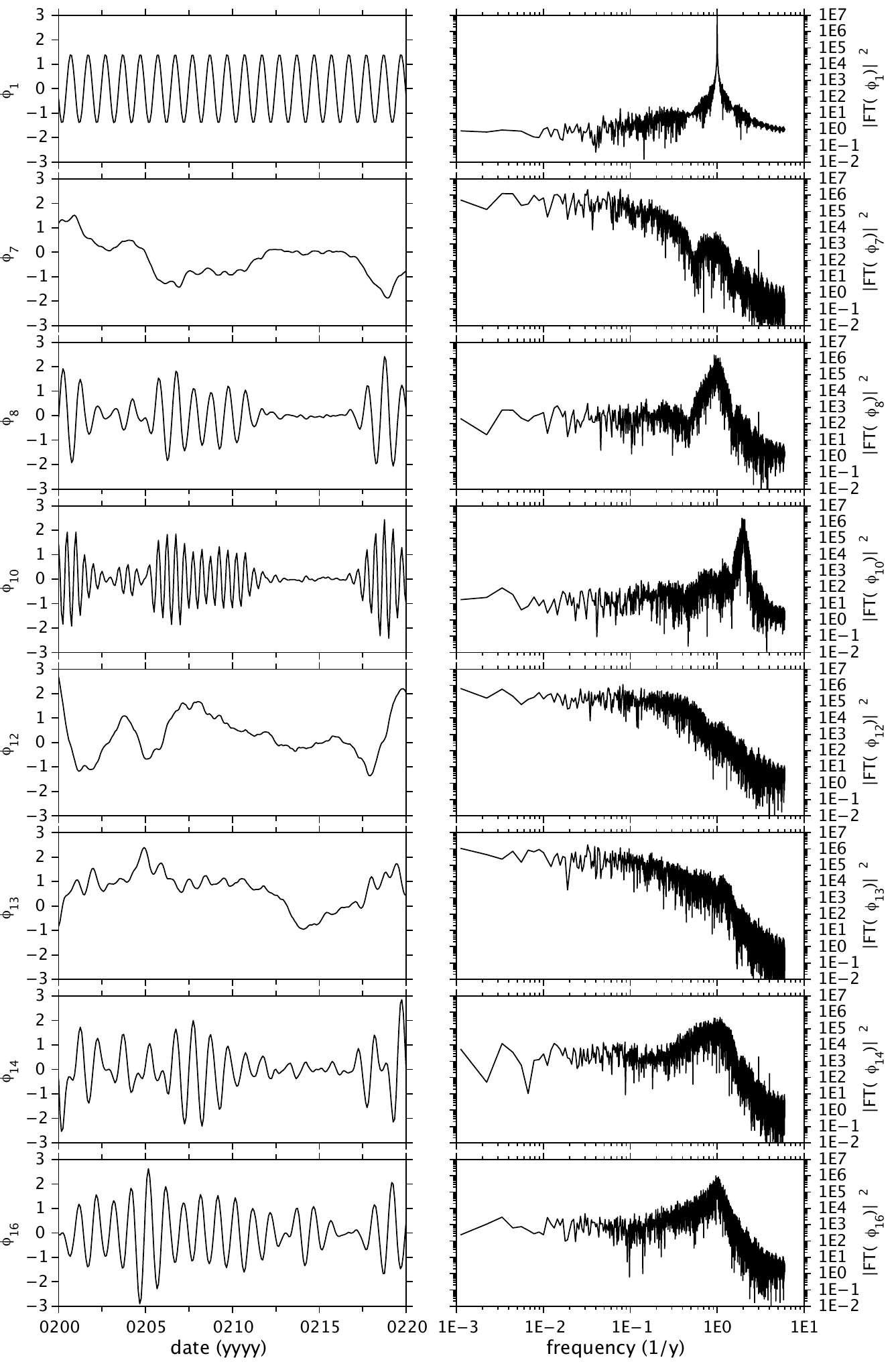}
  \caption{\label{figCCSMPhiCone}Same as Figure~\ref{figCCSMPhiAT} but for the cone kernel with $ \zeta = 0.995 $. The displayed eigenfunctions are (a) annual periodic; (b) PDO; (c,d) annual and semiannual intermittent patterns associated with the PDO; (e) ENSO; (f) NPGO; (g) annual intermittent patterns associated with ENSO; (h) annual intermittent patterns associated with the NPGO.} 
\end{figure}

\begin{figure}
  \centering\includegraphics[scale=.84]{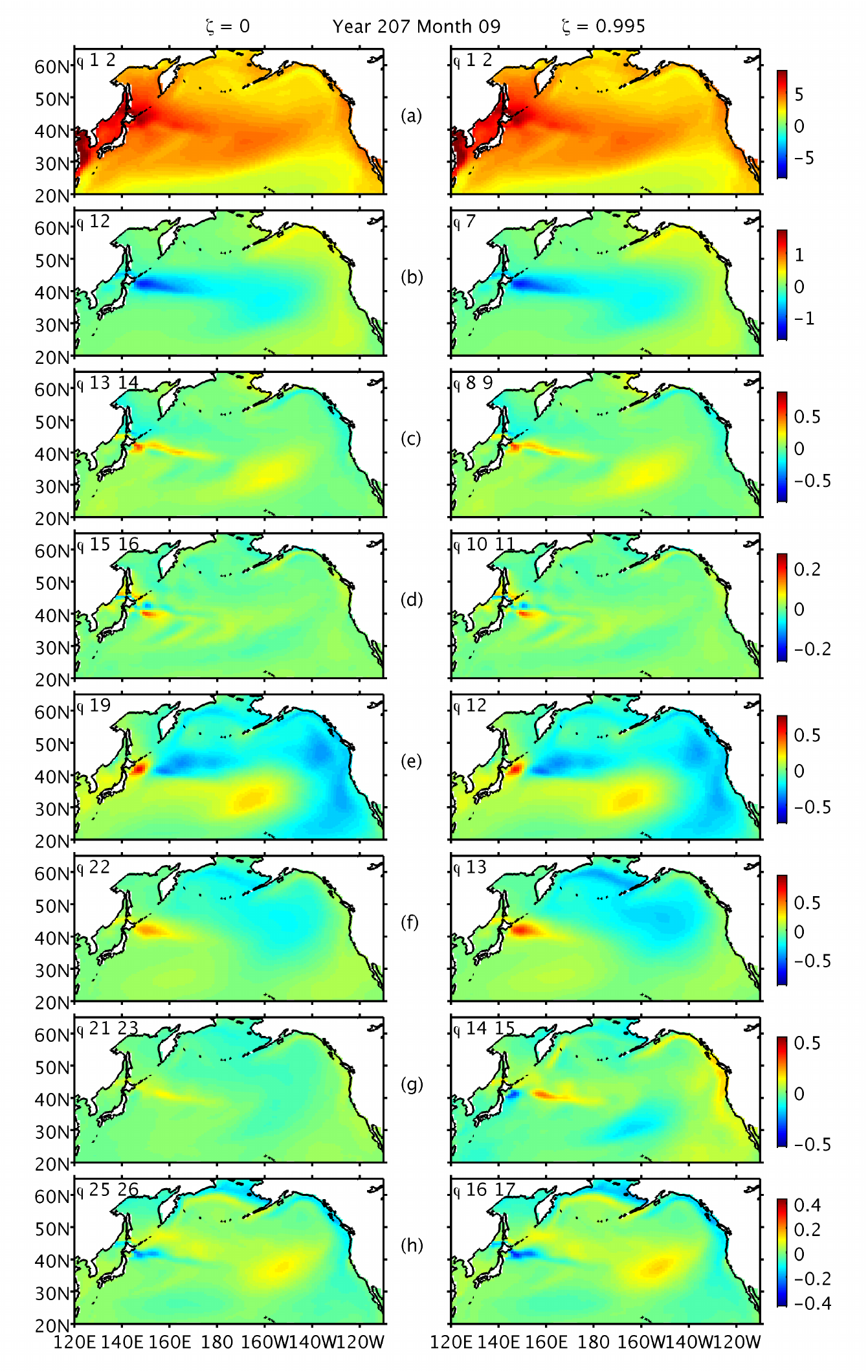}
  \caption{\label{figCCSMSnapshot}Snapshots of spatiotemporal patterns of SST anomalies (in K) for September of simulation year 207 obtained from the eigenfunctions in Figures~\ref{figCCSMPhiAT} and~\ref{figCCSMPhiCone}. The displayed patterns are the (a) annual periodic; (b) PDO; (c,d) annual and semiannual intermittent patterns associated with the PDO; (e) ENSO; (f) NPGO; (g) annual intermittent patterns associated with ENSO; (h) annual intermittent patterns associated with the NPGO. Note the higher amplitude in panels (f--h) for the $ \zeta = 0.995 $ cone kernel compared to the $ \zeta =0 $ cases. See Movie~1 for the dynamic evolution of these patterns and the raw data.}
\end{figure}

Representative diffusion eigenvalues $ \lambda_i $ and eigenfunctions $ \phi_i $ obtained via cone kernels with $ \zeta = 0 $ and $ 0.995 $ are displayed in Figures~\ref{figCCSMLambda}--\ref{figCCSMPhiCone}, where the $ \phi_i $ are represented by the corresponding $ \tilde \phi_i $ time series from~\eqref{eqPhiTS}. To visualize subspaces spanned by selected groups of eigenfunctions  $ \phi = ( \phi_{i_1}, \ldots, \phi_{i_l} ) $ as spatiotemporal patterns, we filter the data using the $ s \times l $ matrix $ \phi $ as a convolution filter; i.e., $ X \mapsto X \pi \phi \phi^ T$, where $ \pi = \diag( \pi_1, \ldots, \pi_s ) $ is the diagonal matrix formed by the inner product weights in~\eqref{eqDiscreteHodge}.  Spatiotemporal patterns of this type are shown in Figure~\ref{figCCSMSnapshot} and Movie~1, which is much more revealing. Note that the eigenfunctions computed using a first-order backward FD scheme (not shown here) exhibit minor changes compared to those in Figures~\ref{figCCSMPhiAT} and~\ref{figCCSMPhiCone}. 

\begin{remark}
  { \ \rm In this application, the variations in the phase space velocity velocity norm $ \lVert \xi_i \rVert $ are of order 10\% (cf.\ section~\ref{secTorus}). As a result, the influence of the local scaling by $ \lVert \xi_i \rVert \lVert \xi_j \rVert $ in this case is rather weak, and correspondingly the eigenfunctions obtained via the $ \zeta = 0 $ cone kernel and the isotropic Gaussian kernel (not shown here) do not differ significantly. Here, the main difference between the isotropic and cone kernels is due to the directionality of the dynamical flow, which is influenced strongly by the annual cycle.}
\end{remark} 

\subsection{\label{secSSTFeatures}Spatiotemporal patterns extracted by diffusion eigenfunctions} 
    
In general, the  eigenfunctions can be grouped into three families according to the timescales present in the $ \tilde \phi_i $ time series, which we refer to as periodic, low-frequency, and intermittent \cite{GiannakisMajda12b}. Figures~\ref{figCCSMPhiAT} and~\ref{figCCSMPhiCone} display examples from each family. The periodic eigenfunctions come in doubly-degenerate pairs, and are dominated by a single frequency component which is an integer multiple of the annual cycle. Due to the presence of strong seasonality, one would expect the data manifold to have the structure of a circle along one of its dimensions---the presence of the periodic eigenfunctions is consistent with this picture. The low-frequency family is characterized by red-noise type frequency spectra, featuring significant power at interannual to decadal timescales. These basis functions describe familiar low-frequency patterns of North Pacific SST variability, including the PDO, ENSO, and NPGO [Panels (b), (e), and (f) in Figure~\ref{figCCSMSnapshot} and Movie~1], which are also accessible through SSA algorithms or PCA of seasonally-detrended data. The intermittent eigenfunctions have the structure of amplitude-modulated wavetrains consisting of a periodic carrier signal with integer frequency multiples of the seasonal cycle, which is modulated by a low-frequency envelope. In certain cases (though not always) the low-frequency envelope is correlated strongly with a low-frequency eigenfunction. In Figure~\ref{figCCSMPhiCone}, for instance, the low-frequency envelopes of the intermittent eigenfunctions $ \{ \phi_8, \phi_{10} \} $, $ \phi_{14} $, and $ \phi_{16} $ correlate strongly with the low-frequency eigenfunctions $ \phi_7 $, $ \phi_{12} $, and $ \phi_{13} $, respectively. 

Similarly to the periodic eigenfunctions, the intermittent eigenfunctions arise in near-degenerate pairs (see Figure~\ref{figCCSMLambda}).  The corresponding spatiotemporal patterns [Panels (c), (d), (g), and (h) in Figure~\ref{figCCSMSnapshot} and Movie~1] feature propagating structures in regions of the North Pacific with high variability, including the Kuroshio current, the Bering Sea, and the west coast of North America. These patterns carry approximately two orders of magnitude less variance of the raw data than the prominent low-frequency modes. As a result, they are not accessible to variance-greedy PCA-type algorithms. However, these eigenfunctions are crucial in explaining the lagged temporal correlation structure of the data in so-called anomaly reemergence mechanisms \cite{AlexanderEtAl99,BushukEtAl13}, and have also been found to have high skill as external factors in regression models for SST variability \cite{GiannakisMajda12c}.

\subsection{Influence of the angular term}

The periodic, low-frequency, and intermittent families arise in both of the $ \zeta = 0 $ and 0.995 kernels. As is evident in Figures~\ref{figCCSMPhiAT} and~\ref{figCCSMPhiCone}, there exists a reasonably clear qualitative correspondence between the leading eigenfunctions in each case. Nevertheless, the two spectra differ in two important aspects. 

First, the number of eigenfunctions required to capture the prominent features described in section~\ref{secSSTFeatures} is significantly smaller in the $ \zeta = 0.995 $ case than the $ \zeta = 0 $ example. In particular, it takes 17 $ \zeta = 0.995 $ eigenfunctions to describe the annual and semiannual cycles, the PDO, ENSO, and NPGO low-frequency modes, and the associated intermittent modes, versus 26 $ \zeta = 0 $ eigenfunctions. The extra eigenfunctions contained in the $ \zeta = 0 $ spectrum are generally associated with higher harmonics of the seasonal cycle, which appear to be of limited physical significance. Examples of these harmonics can be seen in the eigenvalue spectrum of Figure~\ref{figCCSMLambda}, where the annual and semiannual periodic eigenfunctions ($ \phi_1 $--$ \phi_4 $) are followed by higher harmonics all the way to the Nyquist limit $ 5 $ y$^{-1} $ associated with the one-month sampling interval. Higher-frequency harmonics are also present in the intermittent eigenfunctions. 

The suppression of the high-frequency harmonics from the $ \zeta = 0.995 $ spectrum is consistent with the observation made in section~\ref{secEigenfunctions} that large gradients of the eigenfunctions along the direction of the dynamical vector field $ v $ lead to large Dirichlet form in~\eqref{eqDirichlet}, and correspondingly large eigenvalues. Because a major component of $ v $  is due to the annual cycle, the higher-harmonics incur a large $ v( \phi_i ) $ penalty, and thus are removed from the leading part of the spectrum. We expect this to be a generic feature in datasets with prominent quasi-periodic behavior. 

A further key difference between the $ \zeta = 0 $ and $ \zeta = 0.995 $ eigenfunctions pertains to timescale separation. In particular, as can be seen in Figure~\ref{figCCSMPhiAT}, the $ \phi_{22} $ eigenfunction for $ \zeta = 0 $ representing the NPGO exhibits a mixture of interannual and annual timescales. Likewise, $ \phi_{21} $ has a mixed intermittent--low-frequency character. On the other hand, the corresponding $ \zeta = 0.995 $ eigenfunctions in Figure~\ref{figCCSMPhiCone} ($ \phi_{13} $ and $ \phi_{14} $, respectively) have a well-separated temporal character, with $ \phi_{13} $ acting as a low-frequency envelope of $ \phi_{16} $. In consequence, the corresponding spatiotemporal patterns have larger amplitude and feature more prominent intermittent coherent structures; see, e.g., the portion of the dynamic evolution around simulation year 207 in Figure~\ref{figCCSMSnapshot} and Movie~1. Again, we attribute this behavior to the superior adaptation of the $ \zeta = 0.995 $ eigenfunctions to the dynamical flow.

\section{\label{secConclusions}Conclusions}

In this work, we have developed a family of kernels for data analysis in dynamical systems, which incorporate empirical information about the dynamical flow in phase space. Compared to canonical isotropic Gaussian kernels, the kernels presented here assign higher affinity to pairs of data samples whose relative displacement vector lies on a small-angle cone with axis parallel to the dynamical vector field $ v $ on the data manifold. The latter is estimated through finite differences of time-ordered samples, i.e., without requiring prior knowledge of the equations of motion. Due to the presence of the angular term, we refer to the new kernels as cone kernels.  

Cone kernels also feature a scaling factor introduced heuristically in earlier work on so-called NLSA algorithms \cite{GiannakisMajda12a,GiannakisMajda13}, whose role is to locally decrease (increase) the rate of decay of the kernel in regions of phase space where the dynamical flow is fast (slow) in the sense of the norm $ \lVert v \rVert_g $ of the dynamical vector field in data space. Moreover, the strength of the angular dependence is controlled by a parameter $ \zeta \in [ 0, 1 ) $, such that $ \zeta = 0 $ and $ \zeta \to 1 $ correspond to zero or maximal influence of the angular term. Thus, cone kernels include the earlier NLSA kernels as the special case $ \zeta = 0 $.  

We evaluated the metric tensor $ h $ induced on the data manifold by cone kernels (Lemma~\ref{lemmaHMetric}) using the asymptotic analysis framework of Berry~\cite{Berry13}, and also studied the associated diffusion operator $ \upDelta_\zeta $. By virtue of the $ \lVert v \rVert_g $-dependent scaling, the induced metric is invariant for all $ \zeta $ under conformal transformations of the original ambient space metric $ g $. Moreover, for $ \zeta > 0 $,  $ h $ contracts local distances between points on the data manifold whose relative displacement is parallel to $ v $, becoming degenerate as $ \zeta $ approaches 1. In that regime, $ \upDelta_\zeta $ becomes along $ v $ \cite{ElworthyEtAl10}, in the sense that the associated codifferential operator asymptotically annihilates all differential 1-forms $ w $ with the property $ w( v ) = 0 $ (Lemma~\ref{lemmaAlong}).

Intuitively, one thinks of $ \upDelta_\zeta $ in the $ \zeta \to 1 $ limit as generating diffusions along the integral curves of the dynamical vector field. Because $ v $ and its integral curves do not depend on the ambient space metric, this feature is intrinsic to the dynamical system under study. More generally, as $ \zeta \to 1 $, the action of $ \upDelta_\zeta $ on functions depends on the ambient space metric only through the ratio $ \mu / \lVert v \rVert_g^m $ for an $ m $-dimensional manifold. The latter is invariant under conformal transformations of $ g $, as well as more general transformations. A further important property arising in the $ \zeta \to 1 $ limit is that the Dirichlet form associated with $ \upDelta_\zeta $ depends on the directional derivative of functions along $  v $, as opposed to the canonical dependence on the gradient operator [see~\eqref{eqDirichlet}]. This property has significant bearing on the structure of the corresponding diffusion eigenfunctions, which are useful in a wide range of dimension reduction, signal processing, and learning problems. In particular, the leading eigenfunctions are expected to be adapted to the dynamical system generating the data, in the sense of varying weakly along the integral curves of $ v $. 

We discussed the utility of cone kernels in a suite of numerical experiments involving nonlinear flows on the 2-torus and North Pacific sea surface temperature (SST) data from a comprehensive climate model. In the torus experiments, we explicitly demonstrated the adaptivity of the diffusion eigenfunctions associated with $ \zeta \approx 1 $ cone kernels to the dynamics (Figures~\ref{figPhiI} and~\ref{figPhiII}), as well as the robustness of those eigenfunctions to non-conformal deformations of the data (Figure~\ref{figPhiIDeformed}). In the North Pacific SST experiments, cone kernels were found to have superior feature extraction capabilities, in the sense of requiring fewer basis functions than their $ \zeta = 0 $ counterparts to describe the salient coherent structures of North Pacific SST variability, while also providing better separation between the timescales associated with the annual solar forcing and low-frequency (interannual) variability of the ocean. We attribute this improvement of skill to the ability of $ \zeta \approx 1 $ cone kernels to take into account changes in the direction of the dynamical flow due to the annual cycle. This feature should be generic in datasets with prominent quasiperiodic behavior. 

There are several open questions generated by this work which lie outside the scope of the present paper. On the theory side, a more detailed understanding of the diffusion eigenfunctions obtained from cone kernels would be desirable; e.g., their properties as embedding coordinates \cite{BerardEtAl94,JonesEtAl08,Portegies13}. Furthermore, it would be useful to explore generalizations of the deterministic framework adopted here to stochastic dynamical systems. Such approaches might involve replacing the ambient-space distances and inner products with suitable  statistical metrics \cite{TalmonCoifman13}, retaining the explicit dependence on aspects related to the drift of the system in phase space. Even in the deterministic dynamical system context, potential shortcomings of cone kernels may arise due to sensitivity to observational noise and/or poor performance of the directional term in high intrinsic dimensions. Such scenarios would warrant modification of the cone kernel formulation put forward here, but we expect the general approach of incorporating empirically accessible information about the dynamics in kernel design to remain fruitful.     

\Appendix
\section{Technical results} 
This Appendix contains the proof of Lemmas~\ref{lemmaFD}, \ref{lemmaHMetric}, and~\ref{lemmaAlong}, as well as derivations of a number of results used in the main text.  Whenever convenient, we use the shorthand notation $ \partial_\mu $ to represent partial differentiation $ \partial / \partial u^\mu $ with respect to a manifold coordinate $ u^\mu $. 

\subsection{\label{appFD}Proof of Lemma~\ref{lemmaFD}}

Let $ ( u^1, \ldots, u^m ) $ be a coordinate system defined in a neighborhood of $ a_i \in \mathcal{ M } $, and $ U_1, \ldots, U_m $ be the corresponding coordinate basis vectors of $ T_{a_i} \mathcal{ M } $. Denote the coordinates of the curve $ \Phi_t a_0 $ by $ u^\mu( t ) $, so that
\begin{displaymath}
  v \rvert_{a_i} = \sum_{\mu=1}^m v^\mu U_\mu, \quad \text{with} \quad v^\mu = \left. \frac{ d u^\mu }{ d t } \right \rvert_{t_i}.
\end{displaymath}
Moreover, fix a  basis $ e_1, \ldots, e_n $ of $ \mathbb{ R }^n $ so that $ X_i = F( a_i ) = \sum_{\nu=1}^n X^\nu e_\nu $. The $ \{ e_\nu \} $ basis can be identified with a basis of $ T_{X_i} \mathbb{ R }^n $ via the canonical isomorphism $ \mathbb{ R }^n \simeq T_{X_i} \mathbb{ R }^n $. With this choice of bases, the components of the derivative map $ D F : T_{a_i} \mathcal{ M } \mapsto T_{X_i} \mathbb{ R }^n $ become $ {DF^\nu}_\mu = \partial X^\nu / \partial u^\mu \rvert_{a_i}. $ We then compute
\begin{align*}
  \xi_i & =  \delta_p X_i = \left. \frac{ d X }{ d t } \right \rvert_{t_i } \, \delta t + O( \delta t^{p+1} ) =  \left. \left( \sum_{\mu=1}^m \frac{ \partial X }{ \partial u^\mu }  \frac{ d u^\mu }{ dt } \right) \right \rvert_{t_i} \, \delta t + O( \delta t^{p+1} ) \\
  &= \sum_{\mu=1}^m \sum_{\nu=1}^n \left. \frac{ \partial X^\nu }{ \partial u^\mu } \right \rvert_{a_i} v^\mu e_\nu \, \delta t + O( \delta t^{p+1} ) = DF \, v  \rvert_{a_i}  \delta t+ O( \delta t^{p+1} ), 
\end{align*}
which gives~\eqref{eqFDV}. Equation~\eqref{eqFDVNorm} follows immediately:
\begin{gather*}
  \lVert v \rVert_{g,a_i}^2 = ( D F \, v, D F \, v ) \rvert_{a_i} = ( \xi_i + O( \delta t^{p+1} ), \xi_i + O( \delta t^{p+1} ) ) / \delta t^2 \\
  \implies \lVert v_i \rVert_g =  \lVert \xi_i \rVert / \delta t+ O( \delta t^{p} ). \qquad \endproof
\end{gather*}

\subsection{\label{appKDerivatives}Derivatives of cone kernels}

Let $ \{ e_1, \ldots, e_n \} $ be a basis of $ \mathbb{ R }^n $ and $ e^{*1}, \ldots, e^{*n} $ its dual basis with $ e^{*\mu}( e_\nu ) = {\delta^\mu}_\nu $. Let also $ c_{\mu\nu } = ( e_\mu, e_\nu ) $ be the corresponding matrix elements of the data space inner product, giving the components of $ \Xi $ in~\eqref{eqXiHat} via
\begin{displaymath}
  \Xi = \sum_{\nu=1}^n \Xi_\nu e^{*\nu}, \quad \Xi_\nu = \sum_{\mu=1}^n c_{\mu\nu} \xi^\mu, \quad \text{with} \quad \xi = \sum_{\mu=1}^n \xi^\mu e_\mu.
\end{displaymath}

To evaluate the derivatives of cone kernels  in~\eqref{eqKDerivatives},  it is convenient to write down $ K_{\delta t, \zeta }( a, \exp_a u ) = \exp( - A( u ) B( u ) ) $ with
\begin{gather*}
    A( u ) = \frac{ 1 }{ \lVert \xi \rVert^2 } \left( \lVert \omega \rVert^2 - \zeta \frac{ ( \xi, \omega )^2 }{ \lVert \xi \rVert^2 } \right), \quad B( u ) = \frac{ \lVert \xi \rVert }{ \lVert \xi' \rVert } \sqrt{ C( u ) }, \\
    C( u ) = \frac{ \lVert \xi' \rVert^2 \lVert \omega \rVert^2 - \zeta ( \xi', \omega )^2 }{ \lVert \xi \rVert^2 \lVert \omega \rVert^2 - \zeta( \xi, \omega )^2 }, \quad \omega = F( a' ) - F( a ),  \quad \xi' = \delta_p F( a' ).
\end{gather*}
The function $ A( u ) $ vanishes at $ u = 0 $, and its leading two derivatives with respect to $ u $ can be computed straightforwardly:
\begin{align*}
  \frac{ \partial A }{ \partial u^\mu } &= \frac{ 2 }{ \lVert \xi \rVert^2 } \sum_{\rho,\sigma=1}^n \left( c_{\rho\sigma} \omega^\rho \frac{ \partial \omega^\sigma }{ \partial u^\mu } - \zeta ( \xi, \omega ) c_{\rho\sigma} \xi^\rho \frac{ \partial \omega^\sigma }{ \partial u^\mu } \right) \\
  \frac{ \partial^2 A }{ \partial u^\mu \, \partial u^\nu } &= \frac{ 2 }{ \lVert \xi \rVert^2 } \sum_{\rho,\sigma=1}^n \left( c_{\rho\sigma} \frac{ \partial \omega^\rho }{ \partial u^\nu } \frac{ \partial \omega^\sigma }{ \partial \omega^\mu } + c_{\rho \sigma} \omega^\rho \frac{ \partial^2 \omega^ \sigma }{ \partial u^\mu \, \partial u^\nu } \right)  \\
  & \quad  - \frac{ \zeta }{ \lVert \xi \rVert^2 } \sum_{\alpha,\beta,\rho,\sigma=1}^n \left(  c_{\alpha\beta} c_{\rho\sigma} \xi^\alpha \xi^\rho \frac{ \partial \omega^\beta }{ \partial u^\nu } \frac{ \partial \omega^\sigma }{ \partial u^\nu } - ( \xi, \omega ) c_{\rho\sigma} \xi^\rho \frac{ \partial^2\omega^\sigma }{ \partial u^\mu \, \partial u^\nu } \right).
\end{align*}
Noting that $ \partial \omega^\nu / \partial u^\mu \rvert_{u=0} = {DF^\nu}_\mu $, these expressions lead to
\begin{displaymath}
  \left. \frac{ \partial A }{ \partial u^\mu } \right \rvert_{u=0} = 0, \quad \left. \frac{ \partial^2 A }{ \partial u^\mu \partial u^\nu } \right \rvert_{u=0} = \frac{ 2 }{ \lVert \xi \rVert^2 } \left( g_{\mu \nu} - \zeta \frac{  \xi^*_\mu \xi^*_\nu }{ \lVert \xi \rVert^2 } \right), 
\end{displaymath}
where 
\begin{equation}
  \label{eqGComponents}g_{\mu\nu} = \sum_{\rho,\sigma=1}^n c_{\rho\sigma} \, {DF^\rho}_\mu \, {DF^\sigma}_\nu 
\end{equation}
are the components of the ambient spaced induced metric $ g $ in~\eqref{eqG}, and $ \xi^*_\mu $ are defined in~\eqref{eqXiHat}. Therefore, assuming that $ B $ and its first two derivatives are all bounded at $ u = 0 $, we have
\begin{equation}
  \label{eqKDerivativesAB}
  \begin{aligned}
    \left. \frac{ \partial K_{\delta t, \zeta } }{ \partial u^\mu } \right \rvert_{u=0} &= \left. - K_{\delta t, \zeta}( a, a ) \left(  \frac{ \partial A }{ \partial u^\mu } B + A \frac{ \partial B }{ \partial u^\mu } \right) \right\rvert_{u=0} = 0, \\ 
    \left. \frac{ \partial^2 K_{\delta t, \zeta } }{ \partial u^\mu u^\nu } \right \rvert_{u=0} &= - \left. \left[ \frac{ \partial K_{\delta t, \zeta } }{ \partial u^\mu }  \left(  \frac{ \partial A }{ \partial u^\mu } B + A \frac{ \partial B }{ \partial u^\mu } \right)  \right] \right \rvert_{u=0} \\
    & \quad - \left. K_{\delta t, \zeta}( a, a ) \left( \frac{ \partial^2 A }{ \partial u^\mu \partial u^\nu } B + \frac{ \partial A }{ \partial u^\mu } \frac{ \partial B }{ \partial u^\nu } + \frac{ \partial A }{ \partial u^\nu } \frac{ \partial B }{ \partial u^\mu } + A \frac{ \partial^2 B }{ \partial u^\mu \partial u^\nu } \right) \right \rvert_{u=0} \\
    & \quad = - \left. \frac{ \partial^2 A }{ \partial u^\mu \partial u^\nu }  B \right \rvert_{u=0}.
  \end{aligned}
\end{equation}  
Thus, in order to obtain~\eqref{eqKDerivatives}, it suffices to determine the value $ B( 0 ) $ and check that the first two derivatives of $ B $ are bounded at $ u = 0 $. Now, because $ \xi $ is a smooth function of $ u $ and $ B( 0 ) = \sqrt{ C( 0 ) } $, any pathological behavior of $ B $ and its derivatives at the origin would be caused by $ C $. 

To confirm that this is not the case, consider the non-degenerate type $ ( 0, 2 ) $ tensor $ Q( u ) $ with components
\begin{displaymath}
  Q( u )_{\mu\nu} = ( \lVert \xi \rVert^2 \delta_{\mu\nu} - \zeta \xi^*_\mu \xi^*_\nu ) \rvert_u,
\end{displaymath}
where $ \xi^*_\mu $ is given by~\eqref{eqXiHat}. Using $ Q( u ) $, we expand $ C( u ) $ as follows: 
\begin{gather*}
  C( u ) = \frac{ \sum_{\mu,\nu=1}^m Q( u )_{\mu\nu} u^\mu u^\nu }{ \sum_{\rho,\sigma=1}^m  Q( 0 )_{\rho\sigma} u^\rho u^\sigma } = 1 + \sum_{\alpha=1}^m u^\alpha R^{(1)}_\alpha + \frac{ 1 }{ 2 } \sum_{\alpha,\beta=1}^m u^\alpha u^\beta R^{(2)}_{\alpha\beta} + O( \lVert u \rVert^3 ), \\
  R^{(1)}_\alpha = \frac{ \sum_{\mu,\nu=1}^m u^\mu u^\nu \partial_\alpha Q_{\mu\nu}( 0 ) }{ \sum_{\rho,\sigma=1}^m u^\rho u^\sigma Q_{\rho\sigma}(0)}, \quad R^{(2)}_{\alpha\beta} = \frac{ \sum_{\mu,\nu=1}^m u^\mu u^\nu \partial_\alpha \partial_\beta Q(0)_{\mu\nu} } { \sum_{\rho,\sigma=1}^m u^\rho u^\sigma Q(0)_{\rho\sigma}}.
\end{gather*}
The tensor components $ R_\alpha^{(1)} $ and $ R^{(2)}_{\alpha\beta} $ are bounded in absolute value through the inequalities
\begin{displaymath}
  \lvert R_\alpha^{(1)} \rvert \leq \lvert\Lambda_\text{max}^{(1)} / \Lambda_\text{min}^{(1)} \rvert, \quad \lvert R_{\alpha\beta}^{(2)} \rvert \leq \lvert\Lambda_\text{max}^{(2)} / \Lambda_\text{min}^{(2)} \rvert,
\end{displaymath}
where $ \Lambda_\text{max}^{(i)} $ ($ \Lambda_\text{min}^{(i)} $) are the largest (smallest) generalized eigenvalues 
\begin{displaymath}
  \partial_\alpha Q(0)_{\mu\nu} y = \Lambda^{(1)} Q(0)_{\mu\nu} y, \quad \partial_\alpha \partial_\beta Q(0)_{\mu\nu} y = \Lambda^{(2)} Q(0)_{\mu\nu} y  
\end{displaymath}
in absolute value. Thus, $ C( 0 ) = B( 0 ) = 1 $, and the leading two derivatives of $ C $ and $ B $ are all bounded, leading through~\eqref{eqKDerivativesAB} to~\eqref{eqKDerivatives}.  

\subsection{\label{appHMetric}Proof of Lemma~\ref{lemmaHMetric}} 

We prove this Lemma following the analysis of Berry~\cite{Berry13} connecting the induced metric to the Hessian of the kernel. 

First, fix a basis $ U_1, \ldots, U_m $ of $ T_a \mathcal{ M } $, and consider the corresponding exponential coordinates $ u^1, \ldots, u^m $ with respect to the ambient-space induced metric $ g $ in~\eqref{eqG}. For small-enough sampling interval $ \delta t $, the $ m $-dimensional Euclidean ball $ B_{\delta t } $ of radius $ \lVert v \rVert_{g,a} \, \delta t $ centered at the origin maps diffeomorphically to a neighborhood $ \mathcal{ B }_{\delta t} = \exp_a B_{\delta t } \subset \mathcal{ M } $ of $ a $. In these coordinates, the volume element $ d \mu $ associated with $ g $ becomes $  \mu \rvert_{ \exp_a a } = \det D \exp_a \rvert_u \, dU^{*1} \wedge \cdots \wedge dU^{*m} $, where $ dU^{*\mu} $ are the dual basis vectors to $ U_\mu $, and $ D \exp_a \rvert_u : T_u T_a \mathcal{ M } \mapsto T_a \mathcal{ M } $ is the derivative of the exponential map,
\begin{displaymath}
  D \exp_a \rvert_u( u' ) = \lim_{\tau \to 0} ( \exp_a( \tau u' + u ) - \exp_a( u ) ) / \tau.
\end{displaymath}
 Moreover, because the cone kernel $ K_{\delta t, \zeta }( a, \exp_a u ) $ in~\eqref{eqKCone} decays exponentially away from $ u = 0 $ for $ \zeta \in [ 0, 1 ) $, it is possible to restrict the integral over $ \mathcal{ M } $ in~\eqref{eqHIntegral} to an integral over $ B_{\delta t } $, incurring an exponentially small error  as $ \delta t \to 0 $ (e.g., Lemma~1 in \cite{BelkinNiyogi08}). For our purposes, it suffices to consider terms up to order $ \delta t^{3/2} $ in the asymptotics \cite{CoifmanLafon06,Berry13}, i.e.,
\begin{equation}
  \label{eqHExp}
  H_{\delta t, \zeta}f( a ) = \frac{ 1 }{ \delta t^m } \int_{B_{\delta t}} K_{\delta t, \zeta} ( a, \exp_a u ) \tilde f( u )  \det D \exp_a \rvert_u \, du^1 \cdots du^m + O( \delta t^{3/2} ),
\end{equation}
where $ \tilde f( u ) = f( \exp_a u ) $. Taylor-expanding this expression about $ u = 0 $ using~\eqref{eqKDerivatives} and~\eqref{eqKHessianV} with $ p \geq 4 $, we obtain 
\begin{equation}
  \label{eqHB}
  H_{\delta t, \zeta}f( a ) = \frac{ 1 }{ \delta t^m } \int_{B_{\delta t}} \left( 1 - \frac{ h( u, u ) }{ \delta t^2 } \right) \det D \exp_a \rvert_u \, du^1 \cdots du^m + O( \delta t^{3/2} ),
\end{equation}
where $ h $ is the Riemannian metric specified in~\eqref{eqHMetric}. Note that $  h = \sum_{\mu,\nu=1}^m h_{\mu\nu} U^{*\mu} U^{*\nu} $, where 
\begin{equation}
  \label{eqHMetricComponents}
  h_{\mu\nu} = \frac{ 1 }{ \lVert v \rVert_g^2 } \left( g_{\mu\nu} - \zeta \frac{ v^*_\mu v^*_\nu }{ \lVert v \rVert_g^2 } \right), \quad v^*_\mu = \sum_{\nu=1}^m g_{\mu\nu} v^\mu,
\end{equation}
and $ g_{\mu\nu } $ are the components of $ g $ in~\eqref{eqGComponents}.

Next, let $ \widehat{\exp}_a $  and $ \hat u^\mu $ be the exponential map and the corresponding normal coordinates associated with $ h $, respectively. For small-enough $ \delta t $, there exists $ \hat B_{\delta t} \subset \mathbb{ R }^m $ such that $ \widehat{ \exp }_a \hat B_{\delta t} = \mathcal{ B }_{\delta t} $ and 
\begin{displaymath}
  du^1 \cdots du^m = \det D \exp_a^{-1} \widehat{\exp}_a \rvert_{\hat u} \, d\hat u^1 \cdots d\hat u^m.
\end{displaymath}
Equation~\eqref{eqHB} can then be expressed as
\begin{equation}
  \label{eqHBHat}
  H_{\delta t, \zeta} f( a ) = \frac{ 1 }{ \delta t^m} \int_{\hat B_{\delta t}} \left( 1 - \frac{ \lVert \hat u \rVert^2 }{ \delta t^2 } \right) \, \hat f( \hat u ) \det D \widehat{\exp}_a \rvert_{\hat u} \, d\hat u^1 \cdots d\hat u^m + O( \delta t^{3/2}),
\end{equation}
with $ \lVert \hat u \rVert^2 = \sum_{\mu=1}^m (\hat u^\mu)^2 $ and $ \hat f( \hat u ) = f( \widehat{\exp}_a \hat u ) $.  Following similar arguments as those used to derive~\eqref{eqHExp}, it is possible to convert~\eqref{eqHBHat} to an isotropic Gaussian integral over $ \hat B_{\delta t } $, i.e.,
\begin{equation}
  \label{eqHBHatGauss}
  H_{\delta t, \zeta} f( a ) = \frac{ 1 }{ \delta t^m } \int_{\hat B_{\delta t} } e^{-\lVert \hat u \rVert^2 / \delta t^2 } \hat f( \hat u ) \det D \widehat{\exp}_a \rvert_{\hat u } \, d\hat u^1 \cdots d\hat u^m + O( \delta t^3 ).
\end{equation}
It then follows from Lemma~8 of~\cite{CoifmanLafon06} (or Lemma~9 of~\cite{BelkinNiyogi08}) that
\begin{align}
  \nonumber H_{\delta t, \zeta} f( a ) &= C_0 f( a ) + \hat C_2 \, \upDelta_{\mathbb{ R }^m} \left( \hat f( 0 )  \det D \widehat \exp_a \rvert_{\hat u = 0} \right) \, \delta t^2 + O( \delta t^3 ) \\
  \label{eqHBHatLapl} &= C_0 ( f( a ) + ( C_2 \kappa( a ) f( a ) + \upDelta_{\zeta} f( a ) ) \, \delta t^2  + O( \delta t^3 ).
\end{align}
Here, $ C_0 $, $  C_2 $, and $ \hat C_2 $ are constants, $ \upDelta_{\mathbb{ R}^m} $ the canonical Laplacian in $ \mathbb{ R }^m $ with $ \upDelta_{\mathbb{ R }^m} \hat f  = \sum_{\mu=1}^m \partial^2 \hat f/ \partial {u^\mu}^2 $, and $ \kappa( a ) $ a function proportional to the scalar curvature of the $ h $ metric at $ a $.  

Inserting~\eqref{eqHBHatLapl} into~\eqref{eqLContinuous}, we then have that asymptotically as $ \delta t \to 0 $ $ \mathcal{ L }_{\delta t, \zeta } f( a )= C \upDelta_{\zeta} f( a ) + O( \delta t ) $ for some constant $ C $ (Theorem~2.4.1 in \cite{Berry13}), which implies the result in~\eqref{eqHMetric}.

To verify~\eqref{eqVolH}, it suffices to compute the determinants of the matrices formed by the $ g_{\mu\nu} $ and $ h_{\mu\nu} $ components and the vector field norm $ \lVert v \rVert_g $  in any convenient coordinate system. In particular, \eqref{eqVolH} follows immediately by evaluating the expressions in~\eqref{eqVNormExp} and~\eqref{eqDetH3} for $ \lVert v \rVert_g $ and $ \det h $ in terms of exponential normal coordinates for $ g $ at $ u = 0 $. $ \qquad \endproof$ 

\subsection{\label{appAlong}Proof of Lemma~\ref{lemmaAlong}}

We work in Riemannian normal coordinates $ u^1, \ldots, u^m $ associated with an orthonormal basis $ U_1, \ldots, U_m $ of $ T_a \mathcal{ M } $ such that $ U_1 $ is aligned with the dynamical vector field $ v $ at $ a \in \mathcal{ M } $. Thus, in this coordinate system we have
\begin{equation}
  \label{eqVExp}
  \begin{aligned}
    v^1 &= v^1( 0 ) + \sum_{\mu=1}^m \partial_\mu v^1( 0 )  u^\mu + O( \lvert u \rVert^2 ), \\
    v^\nu &= \sum_{\mu=1}^m \partial_\mu v^\nu( 0 )  u^\mu + O( \lVert u \rVert^2 ), \quad \nu > 1, \\
    g_{\mu\nu} &= \delta_{\mu\nu} + O( \lVert u \rVert^2 ).
  \end{aligned}
\end{equation}
Let $ w = \sum_{\mu=1}^m w_\mu U^{*\mu} $ be a $ C^1 $ 1-form field with components $ w_\mu $ in the dual basis $ U^{*\mu} $ to $ U_\mu $. Our objective is to evaluate $ \delta_\zeta w $ from~\eqref{eqLapl} at $ a $ using the standard expression for the codifferential in  a local coordinate system (e.g., \cite{Rosenberg97}),
\begin{align}
  \nonumber \delta_\zeta w  &=  - \left. \frac{ 1 }{ \sqrt{\det h }} \sum_{\mu,\nu=1}^m \partial_\mu( h^{-1,\mu\nu} \sqrt{ \det h } w_\nu ) \right\rvert_{u=0} \\
  \label{eqCodiff}
  &= - \sum_{\mu,\nu=1}^m \left. \left( h^{-1,\mu\nu}  \partial_\mu w_\nu + w_\nu \partial_\mu h^{-1,\mu\nu} + \frac{ h^{-1,\mu\nu} w_\nu }{ 2 \det h }  \partial_\mu \det h \right) \right\rvert_{u=0},
\end{align}  
and then collect the dominant terms in $ ( 1 - \zeta ) $. 

First, consider the dual vector field $ v^* = g( v, \cdot ) = \sum_{\mu=1}^m v^*_\mu U^{*\mu } $ defined in~\eqref{eqKHessianV} and its norm $ \lVert v \rVert_g $ with respect to the data space metric $ g $. It follows from~\eqref{eqVExp} that
\begin{align} 
\nonumber  v^*_1 &= v^1( 0 ) + \sum_{\mu=1}^m \partial_\mu v^1( 0 ) u^\mu + O( \lVert u \rVert^2 ), \\
\nonumber  v^*_\nu &= \sum_{\mu=1}^m \partial_\mu v^\nu( 0 ) u^\mu + O( \lVert u \rVert^2 ), \quad \nu > 1 \\
\label{eqVNormExp}  \lVert v \rVert_g &= v^1( 0 ) + \sum_{\mu=1}^m \partial_\mu v^1( 0 )  u^\mu + O( \lVert u \rVert^2 ).
\end{align}
Moreover,
\begin{align*}
  (v^*_1)^2/\lVert v \rVert_g^2 &= 1 + O( \lVert u \rVert^2 ), \\
  v^*_1 v^*_\mu / \lVert v \rVert_g^2 &= \frac{ 1 }{ v^1( 0 ) } \sum_{\rho=1}^m \partial_\rho v^\mu( 0 )  u^\nu + O( \lVert u \rVert^2 ), \quad \mu > 1, \\
  v^*_\mu v^*_\nu / \lVert v \rVert_g^2 &= O( \lVert u \rVert^2 ), \quad \mu,\nu>1.
\end{align*}
With these results, the metric components $ h_{\mu\nu} $ in~\eqref{eqHMetricComponents} can be expressed as
\begin{equation}
  \label{eqHDecomp}
  h_{\mu\nu} = \bar h_{\mu\nu} + h'_{\mu\nu} + O( \lVert u \rVert^2 ),
\end{equation}
with 
\begin{gather*}
  \bar h_{11} = \frac{ 1 - \zeta }{ ( v^1(0) )^2 }, \quad h'_{11} = -\frac{ 2( 1 - \zeta ) }{ ( v^1(0 )^3 ) } \sum_{\rho=1}^\nu \partial_\rho v^1( 0 )  u^\rho, \\
  \bar h_{1\nu} = 0, \quad h'_{1\nu} = - \frac{ \zeta }{ (v^1(0))^3 } \sum_{\rho=1}^m \partial_\rho v^\nu( 0 )  u^\rho, \quad \nu > 1, \\
  \bar h_{\mu\nu} = \frac{ \delta_{\mu\nu } }{ ( v^1( 0 ) )^2 }, \quad h'_{\mu\nu} = - \frac{ 2 \delta_{\mu\nu} }{ ( v^1( 0 ) )^3 } \sum_{\rho=1}^m \partial_\rho v^1( 0 )  u^\rho, \quad \mu,\nu>1.
\end{gather*}

Next, consider the inverse metric $ h^{-1,\mu\nu} $. Asymptotically, we have
\begin{equation}
  \label{eqHInv}
  h^{-1,\mu\nu} = \bar h^{-1, \mu \nu } - h^{\prime\mu\nu} + O( \lVert u \rVert^2 ), \quad \text{with} \quad h^{\prime \mu\nu} = \sum_{\rho,\sigma=1}^m \bar h^{-1,\mu\rho} h'_{\rho\sigma} \bar h^{-1,\sigma \nu}.
\end{equation}
The tensor components appearing in~\eqref{eqHInv} are
\begin{gather*}
  \bar h^{-1,11} = \frac{ (v^1( 0 ))^2 }{ 1 - \zeta }, \quad h^{\prime 1 1} = - \frac{ 2 v^1( 0 ) }{ 1 - \zeta } \sum_{\rho=1}^m \partial_\rho v^1( 0 )  u^\rho, \\
  \bar h^{-1,1\nu} = 0, \quad h^{\prime, 1\nu} = - \frac{ \zeta v^1( 0 ) }{ 1 - \zeta } \sum_{\rho=1}^m \partial_\rho v^\nu( 0 )   u^ \rho, \quad \nu > 1, \\
  \bar h^{-1,\mu\nu} = ( v^1( 0 ) )^2 \delta_{\mu\nu}, \quad h^{\prime \mu\nu} = - 2 v^1( 0 ) \sum_{\rho=1}^m \partial_\rho v^1( 0 )  u^\rho, \quad \mu,\nu > 1.
\end{gather*}
We therefore obtain
\begin{equation}
  \label{eqCodiff1}
  - \left. \sum_{\mu,\nu=1}^m h^{-1,\mu\nu} \partial_\mu w_\nu \right \rvert_{u=0} = - \sum_{\mu,\nu=1}^m \bar h^{-1,\mu \nu} \partial_\mu w_\nu( 0 ) = - \frac{ ( v^1( 0 ) )^2 }{ 1 - \zeta } \partial_1 w_1( 0 ) - ( v^1( 0 ) )^2 \sum_{\nu=1}^m \partial_1 w_\nu( 0 )
\end{equation}
and
\begin{align}
 \nonumber - \left. \sum_{\mu,\nu=1}^m \partial_\mu h^{-1,\mu\nu} w_\nu \right \rvert_{u=0} &= \sum_{\mu,\nu=1}^m \partial_\mu h^{\prime \mu\nu}( 0 ) w_\nu( 0 ) \\
  \nonumber &= - \frac{ 2 v^1( 0 ) }{ 1 - \zeta } \partial_1 v^1( 0 ) w_1( 0 ) - \frac{ \zeta v^1( 0 ) }{ 1 - \zeta } \sum_{\nu=2}^m \left( \partial_1 v^\nu( 0 ) w_\nu( 0 ) + \partial_\nu v^\nu( 0 ) w_1( 0 ) \right) \\ 
  \label{eqCodiff2} & \quad - 2 v^1( 0 ) \sum_{\nu=2}^m \partial_\nu v^1( 0 ) w_\nu( 0 ).
\end{align}

What remains is to compute the contribution to the codifferential from the derivatives of $ \det h $. To carry this calculation, we begin from the standard formula for the determinant, 
\begin{equation}
  \label{eqDetH}\det h = \sum_{\alpha_1, \alpha_2, \ldots, \alpha_m=1}^m \epsilon_{\alpha_1 \alpha_2 \cdots \alpha_m} h_{1\alpha_1} h_{2\alpha_2}\cdots h_{m\alpha_m},
\end{equation}
where $ \epsilon_{\alpha_1 \alpha_2 \cdots \alpha_m} $ is equal to 1 $(-1) $ if $ ( \alpha_1, \ldots, \alpha_m ) $ is an even (odd) permutation of $ ( 1, \ldots, m ) $, and vanishes  if the $ \alpha_i $ are not distinct.  Because the $ O( \lVert u \rVert ) $ off-diagonal components $ h_{\mu\nu} $ occur for $ \mu = 1 $ or $ \nu = 1 $ only, the $ h_{i \alpha_i} $ factors in~\eqref{eqDetH} are $ O( \lVert u \rVert ) $ only for $ \alpha_i \in \{ 1, \alpha_i \} $. Therefore,
\begin{align*}
  \det h &= \sum_{\alpha_1=1}^m \sum_{\alpha_2 \in \{ 1, 2 \} } \sum_{\alpha_3 \in \{ 1, 3 \} } \cdots \sum_{\alpha_m \in \{ 1, m \} } \epsilon_{\alpha_1 \alpha_2 \cdots \alpha_m }  h_{1 \alpha_1} h_{2\alpha_2} \cdots h_{m\alpha_m} + O( \lVert u \rVert^2 ) \\
  & = \sum_{\alpha_1=1}^m \sum_{\alpha_3 \in \{ 1, 3 \} } \cdots \sum_{\alpha_m \in \{ 1, m \} } 
  ( \epsilon_{\alpha_1 1 \alpha_3 \cdots \alpha_m} h_{1 \alpha_1 } h_{ 21 } h_{3 \alpha_3} \cdots h_{m \alpha_m} \\
  & \qquad + \epsilon_{\alpha_1 2 \alpha_3 \cdots \alpha_m } h_{1\alpha_1} h_{22} h_{3 \alpha_3 } \cdots h_{m \alpha_m} ) + O( \lVert u \rVert^2 ) \\
  & = \sum_{\alpha_1=1}^m \epsilon_{\alpha_1 2 3 \cdots m} h_{1\alpha_1} h_{21} h_{33} \cdots h_{mm} \\
  & \qquad +  \sum_{\alpha_1=1}^m \sum_{\alpha_3 \in \{ 1, 3 \} } \cdots \sum_{\alpha_m \in \{ 1, m \} } \epsilon_{\alpha_1 2 \alpha_3 \cdots \alpha_m} h_{1\alpha_1} h_{22} h_{3\alpha_3} \cdots h_{m\alpha m} + O( \lVert u \rVert^2 )
\\ & = - h_{12} h_{21} h_{33} \cdots h_{mm} +  \sum_{\alpha_1=1}^m \sum_{\alpha_3 \in \{ 1, 3 \} } \cdots \sum_{\alpha_m \in \{ 1, m \} } \epsilon_{\alpha_1 2 \alpha_3 \cdots \alpha_m} h_{1\alpha_1} h_{22} h_{3\alpha_3} \cdots h_{m\alpha_m} \\
& \qquad + O( \lVert u \rVert^2 ).
\end{align*} 
Repeating this decomposition $ m -2 $ times, we obtain 
\begin{equation}
  \label{eqDetH2}
  \det h = h_{11} h_{22} \cdots h_{mm} + ( - h_{12} h_{21} h_{33} \cdots h_{mm} + h_{13} h_{22} h_{31} h_{44} \cdots h_{mm} + \cdots ) + O( \lVert u \rVert^2 ).
\end{equation}
By virtue of~\eqref{eqHDecomp}, the first term in the right-hand side of~\eqref{eqDetH2} is $ O( \lVert u \rVert ) $, whereas each of the terms in the parentheses is $ O( \lVert u \rVert^2 ) $. Specifically,
\begin{equation}
  \label{eqDetH3}
  \det h = \frac{ 1 - \zeta }{ ( v^1( 0 ) )^{2m} } \left( 1 - \frac{ 2 }{ v^1( 0 ) }\sum_{\rho=1}^m \partial_\rho v^1( 0 )  u^\rho \right)^m + O( \lVert u \rVert^2 ),
\end{equation}
which implies that
\begin{displaymath}
  \left. \frac{ \partial_\mu \det h }{ \det h } \right \rvert_{u=0} = - \frac{ 2 m }{ v^1( 0 ) }  \partial_\mu v^1( 0 )
\end{displaymath}
and
\begin{equation}
  \label{eqCodiff3}
  - \left. \frac{ h^{-1,\mu\nu} w_\nu }{ 2 \det h } \partial_\mu \det h \right\rvert_{u=0} =  \frac{ m v^1( 0 ) }{ 1 - \zeta }  \partial_1 v^1( 0 )  w_1( 0 ) + m v^1( 0 ) \sum_{\nu=2}^m \partial_\nu v^1( 0 ) w_\nu( 0 ). 
\end{equation}

Using~\eqref{eqCodiff1}, \eqref{eqCodiff2}, and~\eqref{eqCodiff3}, it is possible to write down an expression for the codifferential $ \delta_\zeta w  $ in the $ u^\mu $ coordinate basis. For the purpose of proving Lemma~\ref{lemmaAlong}, it suffices to consider the dominant, $ O( ( 1 - \zeta )^{-1} ) $, terms, namely
\begin{align}
  \nonumber
  \delta_\zeta( w ) &= - \frac{ ( v^1( 0 ) )^2 }{ 1 - \zeta }  \partial_1 w_1( 0 ) - \frac{ v^1( 0 ) }{ 1 - \zeta } \sum_{\nu=1}^m w_\nu( 0 )  \partial_1 v^\nu( 0 ) - \frac{ v^1( 0 ) }{ 1 - \zeta } \sum_{\nu=1}^m \partial_\nu v^\nu( 0 )  w_1( 0 ) \\ 
   \label{eqCodiff4}  &+ \frac{ m v^1( 0 ) }{ 1 - \zeta }  \partial_1 v^1( 0 )  w_1( 0 ) + O( ( 1 - \zeta )^0 ).
\end{align}
Noting the relations
\begin{gather*}
  w( v ) \rvert_{u=0} = \sum_{\mu=1}^m ( v^\mu w_\mu ) \rvert_{u=0} = w_1( 0 ) v^1( 0 ),\\
  \begin{aligned}
    v( w( v ) ) \rvert_{u=0} &= \sum_{\mu,\nu=1}^m v^\mu \partial_\mu( w_\nu v^\nu ) \rvert_{u=0} = \sum_{\mu,\nu=1}^m ( v^\mu v^\nu \partial_\mu w_\nu + v^\mu w_\nu \partial_\mu v^\nu ) \rvert_{u=0} \\
    &= ( v^1( 0 ) )^2 \partial_1 w_1( 0 ) + v^1( 0 ) \sum_{\nu=1}^m v^1( 0 ) w_\nu( 0 )  \partial_1 v^\nu( 0 ),
  \end{aligned}\\
    \divr_\nu v \rvert_{u=0} = \left. \frac{ 1 }{ \sqrt{ \det h } } \sum_{\nu=1}^m \partial_\nu ( \sqrt{ \det h } v^\mu ) \right\rvert_{u=0} = \sum_{\nu=1}^m \partial_\nu v^\nu( 0 ) - m \partial_1 v^1( 0 ), \\
  \divr_\nu[ w( v ) v ] = w( v ) \divr_\nu v + v( w( v ) ),
\end{gather*}
%v( \lVert v \rVert_g ) = \sum_{\mu,\nu,\rho=1}^m \left( v^\mu  \partial_\mu( v^\nu v^\rho g_{\nu\rho} )^{1/2} \right)_{u=0} = \sum_{\mu,\nu,\rho=1}^m \lVert v \rVert_g^{-1} v^\mu \partial_\mu v^\nu  v^\rho g_{\nu\rho} \rvert_{u=0} = v^1( 0 )  \partial_1 v^1( 0 ),
the asymptotic expansion~\eqref{eqCodiff4} can be put in the covariant (basis-independent) form~\eqref{eqAlong}, thus proving the Lemma. $ \qquad \endproof $  

\subsection*{Acknowledgments} The author would like to thank T.\ Berry, M.\ Bushuk, J.\ Harlim, and A.\ Majda for stimulating discussions. This work was supported by ONR DRI grant N00014-14-1-0150 and ONR MURI grant 25-74200-F7112.

\end{document}